\documentclass[12pt, a4paper]{amsart}
\usepackage[dvips]{epsfig}
\usepackage{amsmath}
\usepackage{amssymb}
\usepackage{amsfonts}
\usepackage{bbm}
\usepackage{amsthm}
\usepackage{amsbsy}
\usepackage{amsgen}
\usepackage{amscd}
\usepackage{amsopn}
\usepackage{amstext}
\usepackage{amsxtra}
\usepackage[utf8]{inputenc}
\usepackage{enumerate}
\usepackage{hyperref}
\usepackage{lineno}
\usepackage{marginnote}
\usepackage{verbatim}
\usepackage{calrsfs}
\usepackage{times, graphicx}
\usepackage{layout}
\usepackage{subcaption}
\usepackage{eso-pic}
\usepackage{lastpage}
\usepackage{dsfont}

\DeclareMathAlphabet{\pazocal}{OMS}{zplm}{m}{n}
\usepackage[foot]{amsaddr}

  \evensidemargin
\oddsidemargin
\numberwithin{equation}{section}

\makeatletter
\newcommand*\rel@kern[1]{\kern#1\dimexpr\macc@kerna}
\newcommand*\widebar[1]{%
  \begingroup
  \def\mathaccent##1##2{%
    \rel@kern{0.8}%
    \overline{\rel@kern{-0.8}\macc@nucleus\rel@kern{0.2}}%
    \rel@kern{-0.2}%
  }%
  \macc@depth\@ne
  \let\math@bgroup\@empty \let\math@egroup\macc@set@skewchar
  \mathsurround\z@ \frozen@everymath{\mathgroup\macc@group\relax}%
  \macc@set@skewchar\relax
  \let\mathaccentV\macc@nested@a
  \macc@nested@a\relax111{#1}%
  \endgroup
}
\makeatother

\theoremstyle{definition}

\newcommand{\vol}{\operatorname{Vol}}

\newcommand{\diam}{\operatorname{diam}}
\newcommand{\Dc}{\mathcal{D}}

\newcommand{\Gc}{\mathcal{G}}

\newcommand{\Sc}{\pazocal{S}}
\newcommand{\Nc}{\pazocal{N}}

\newcommand{\Rc}{\mathcal{R}}

\newcommand{\Qc}{\pazocal{Q}}

\newcommand{\Hc}{\mathcal{H}}

\newcommand{\proj}{\operatorname{Proj}}

\newcommand{\meas}{\operatorname{meas}}

\newcommand {\E} {\mathbb{E}}
\newcommand {\M} {\pazocal{M}}
\newcommand {\R} {\mathbb{R}}

\newcommand {\Z} {\mathbb{Z}}

\newcommand {\Cc} {\pazocal{C}}
\newcommand {\Lc} {\mathcal{L}}

\newcommand {\C} {\mathbb{C}}

\newcommand {\Zc} {\mathcal{Z}}

\newcommand {\Tb} {\mathbb{T}}
\newcommand {\var} {\operatorname{Var}}
\newcommand {\len} {\operatorname{len}}
\newcommand {\cov} {\operatorname{Cov}}
\newcommand {\corr} {\operatorname{Corr}}
\newcommand {\Vc} {\mathcal{V}}

 \allowbreak

\begin{document}

\title[On the nodal structures of random fields]{On the nodal structures of random fields - a decade of results}
\author{Igor Wigman}

\address{Department of Mathematics, King's College London, UK}\email{igor.wigman@kcl.ac.uk}

\begin{abstract}
We survey a decade worth of work pertaining to the nodal structures of random fields, with emphasis on the transformative techniques that shaped the field.
\end{abstract}

\date{\today}
\maketitle

\section{Introduction}

In $2011$ I published the survey ~\cite{W Dartmouth} of results pertaining to the distribution of local functionals of some Gaussian ensembles of random fields,
a subject that, in hindsight, was only in its early stage. Some $10$ years later, one has, what I believe, a complete understanding of the said subject, and it seems like a good
time for another survey of this research direction, to be written. While this manuscript does not mean to contain an exhaustive list of papers on the subject, it aims to describe the transformative or principal results that shaped the field, and outline the main techniques of proofs. We will outline the essence of the arguments behind the proofs of selected results, and emphasize the subtleties of the problems.

\subsection{Euclidean vs. Riemannian setting}

\label{sec:Eucl Riem}

For $n\ge 1$ let $F:\R^{n}\rightarrow\R$ be a random field, that, unless specified otherwise, will be assumed stationary Gaussian, a.s. $C^{\infty}$-smooth.
One is interested in the geometry of $F$, and, importantly, the structure of the (random) nodal set $\Nc_{F}:=F^{-1}(0)\subseteq \R^{n}$. It is usual to restrict
$F|_{B(R)}$ to a large centred ball $B(R)\subseteq \R^{n}$, and study the geometry of the restricted field in the limit $R\rightarrow\infty$; the
corresponding restricted nodal set is $\Nc_{F}\cap B(R)$. Associated to
$F|_{B(R)}$ one distinguishes a class of quantitative {\em local} properties, that could be described as {\em additive} functionals of the domain.
For example, for $n=1$ the number of {\em zeros} of $F$ is a local property, since if $I_{1}=(a,b)$ and $I_{2}(c,d)$ are two {\em disjoint} intervals, then the number of zeros on $I_{1}\cup I_{2}$ is
the sum of zeros of $F$ on $I_{1}$ and on $I_{2}$. Similarly, for $n=2$ (resp. $n\ge 2$), the nodal length (resp. $(n-1)$-volume) of $F|_{B(R)}$, i.e. the length of the smooth
curve $F^{-1}(0)$ could be thought of as the functional $$\Lc_{F}:\mathcal{D}\mapsto \Lc_{F}\mathcal{D}\in \R_{\ge 0},$$ where $\Lc_{F}\mathcal{D}$ is the nodal length of $F$ restricted to
a domain $\mathcal{D}\subseteq \R^{2}$.

Any quantitative property of the geometry of $F$ that fails to be local is said to be {\em nonlocal}.
Examples of important nonlocal properties of $F$ include the number of {\em nodal components} of $F|_{B(R)}$ (i.e. the connected components of $F^{-1}(0)\cap B(R)$), or the percolation properties of excursion sets $F^{-1}(0,+\infty)$ on rectangles $[a,b]\times [c,d]\subseteq \R^{2}$ and their dilates.
The number of nodal components of $F|_{B(R)}$ satisfies an intermediate property of {\em semi-locality} ~\cite{NS1,NS2}. That is, as $R\rightarrow\infty$, with {\em high probability},
{\em most} nodal components of $F|_{B(R)}$ lie in some radius-$r$ ball, with $r>0$ growing {\em slowly} with $R$. In this survey we will only focus our attention to the local properties of ensembles
of random fields in a number of scenarios, especially their nodal volumes, though we will also mention the {\em critical values}, and {\em nodal intersections} against smooth curves.

A particularly important planar random field is Berry's Random Wave Model (RWM), that is the centred Gaussian random field $u:\R^{2}\rightarrow\R$ of monochromatic isotropic waves,
uniquely prescribed by the covariance function
\begin{equation}
\label{eq:RWM covar}
r_{RWM}(x,y)=r_{RWM}(x-y):= \E[u(x)\cdot u(y)] = J_{0}(|x-y|) ,
\end{equation}
with $J_{0}$ Bessel's J function of order $0$.
The said random field, used as a model for {\em ocean waves} by Longuet-Higgins, is believed ~\cite{Berry1977} to represent the (deterministic) Laplace eigenfunctions on generic chaotic surfaces,
in the high energy limit.

Equivalently, one may define $u$ via its {\em spectral measure}: It is the arc length measure of the unit circle $\Sc^{1}\subseteq\R^{2}$, easily extended to higher dimensions.
In general, the spectral measure $\rho=\rho_{F}$ of a {\em stationary} random field $F:\R^{n}\rightarrow\R$ is the Fourier transform of its covariance function $$r_{F}(x,y)=r_{F}(y-x):=\E[F(x)\cdot F(y)],$$
$x,y\in\R^{n}$, thought of as a function of a single variable $r_{F}(z)$. That $\rho$ is indeed a probability measure is the statement of Bochner's theorem. The spectral measure of $F$
uniquely prescribes its law, and oftentimes, rather than imposing some constraints on the covariance function $r_{F}$, one imposes conditions on the law of $F$ via $\rho_{F}$. For example, the support
of $\rho_{F}$ having a point in its interior is a strong non-degeneracy condition on the law of $F$, implying, in particular, that the distribution of $(F(0),\nabla F(0))\in\R^{n+1}$ is non-degenerate,
and, by the stationarity of $F$, so is the distribution of $(F(x),\nabla F(x))\in\R^{n+1}$ for every $x\in\R^{n}$.

Other than the asymptotic geometry of a fixed random field over expanding subdomains of $\R^{n}$, one is interested
in {\em ensembles} of Gaussian random fields defined on a Riemannian manifold $(\M,g)$: $f_{n}:\M\rightarrow\R$, where $n$ belongs to a discrete set of indexes.
In many interesting cases, $f_{n}$ admits scaling limits, defined via the tangent plane, with stationary or even isotropic limit Euclidean random field.
One in particular important such ensemble with a scaling random field is that of {\em band-limited} functions (see e.g. ~\cite{SarWig}) on
a smooth compact Riemannian manifold $\M$, defined in the following section.

\subsection{Random band-limited functions}
\label{sec:band-lim functions}

Let $\Delta$ be the Laplace-Beltrami operator on $\M$ with Dirichlet boundary conditions (say), known to have a purely discrete spectrum, and denote the corresponding
(negative) Laplace eigenvalues $\{\lambda_{i}^{2}\}_{i\ge 0}$ with associate orthonormal basis of $L^{2}(\M)$ consisting of Laplace eigenfunctions $\{\varphi_{i}\}_{i\ge 0}$, i.e. $$\Delta \varphi_{i}+\lambda_{i}^{2}\varphi_{i}=0.$$ Fix a number $\alpha\in [0,1]$, and for the spectral parameter $\lambda\gg 0$ define the random field of random $\alpha$-band limited functions
\begin{equation}
\label{eq:f band lim def}
f_{\lambda}(x) = f_{\M,\alpha;\lambda}(x) = \frac{1}{\sqrt{V(\lambda)}}\sum\limits_{\alpha\cdot\lambda\le \lambda_{i}\le \lambda}c_{i}\varphi_{i}(x),
\end{equation}
indexed by $x\in\M$, where $c_{j}$ are i.i.d. standard Gaussian, and $V(\lambda)>0$ is a convenience pre-factor chosen so that to make $f_{\lambda}(\cdot)$ asymptotically {\em univariate}. For $\alpha=1$, the range of the
summation on the r.h.s. of \eqref{eq:f band lim def} should be understood as $$\lambda-\eta(\lambda)\le\lambda_{i}\le \lambda,$$ with $\eta(\lambda)=o_{\lambda\rightarrow\infty}(\lambda)$ but\footnote{In fact, $\eta(\lambda)$ does not need to grow, under some more restrictive assumptions on $\M$ of geometric nature, e.g. if the measure of geodesic loop directions through $x\in \M$ is $0$ for a.a. $x\in\M$, see a discussion in ~\cite{CH} and the references therein, in particular ~\cite{Safarov}. The case $\eta\equiv 1$ is the most interesting, but also the most difficult
one from the microlocal analytic point of view.}
$\eta(\lambda)\rightarrow\infty$ - these are the {\em ``monochromatic"} waves on $\M$.

One may identify the covariance function of $f_{\lambda}(\cdot)$, uniquely prescribing the law of $f_{\lambda}$, as the {\em spectral projector} $$r_{\lambda}(x,y)=r_{\M,\alpha;\lambda}(x,y):=
\frac{1}{V(\lambda)}\sum\limits_{\alpha\cdot\lambda\le \lambda_{i}\le \lambda}\varphi_{i}(x)\cdot \varphi_{i}(y).$$
It is well known in the microlocal analysis literature that $r_{\lambda}(x,y)$, and thereby $f_{\lambda}(\cdot)$, scale around every {\em reference point} $x_{0}\in\M$, in the following sense.
Let $R>0$ be a number, and assuming that $R/\lambda$ is less than the injectivity radius of $x_{0}$, the (scaled) exponential map $$\exp_{x_{0}}(\cdot/\lambda):B(R)\rightarrow \M$$ assigning
$z\mapsto \exp_{x_{0}}(z/\lambda) \in\M$ is a {\em bijection} between the (Euclidean) ball $B(R)\subseteq T_{x_{0}}\M\cong \R^{n}$ and a neighbourhood of $x_{0}$ in $\M$. We define the scaled Gaussian random field on $B(R)$:
\begin{equation}
\label{eq:band-lim scale g}
g_{\lambda}(z) := f_{\lambda}(\exp_{x_{0}}(z/\lambda)),
\end{equation}
whose covariance function
\begin{equation*}
\widetilde{r}_{\lambda}(z,w) = \frac{1}{V(\lambda)}\sum\limits_{\alpha\cdot\lambda\le \lambda_{i}\le \lambda}\varphi_{i}(\exp_{x_{0}}(z/\lambda))\cdot \varphi_{i}(\exp_{x_{0}}(w/\lambda))
\end{equation*}
converges, uniformly on $B(R)\times B(R)$, together with any finite number of its mixed derivatives, to the Fourier transform of the characteristic function of the annulus $\{\xi\in\R^{n}:\: \alpha\le \|\xi\|\le 1\}$:
\begin{equation}
\label{eq:covar band-lim->covar infty}
\widetilde{r}_{\lambda}(z,w)\rightarrow \widetilde{r}_{\alpha;\infty}:=\int\limits_{\alpha\le |\xi|\le 1} \exp(2\pi i \langle z-w,\xi\rangle ) d\xi,
\end{equation}
with effective control over the error term in terms of power decay in $\lambda$.
For $\alpha=1$, the said characteristic function of the annulus is understood in the limit sense as the hypersurface volume measure of the unit hypersphere $\|\xi\|=1$.

The above means that, for every $R>0$, as $\lambda\rightarrow\infty$, the (not stationary) random field $g_{\lambda}(\cdot)$ converges to the stationary isotropic Gaussian random field $g_{\alpha;\infty}$ on $B(R)$ whose covariance function is $\widetilde{r}_{\alpha;\infty}$. For $\alpha=1$, the random field $g_{\alpha;\infty}$ is identified as Berry's RWM with $n=2$, and, more generally, for $\alpha=1$, $n\ge 2$ these are the monochromatic waves on $\R^{n}$.
These scaling properties allow for inferring results on $f_{\lambda}$ from the corresponding results on $g_{\infty}$, in Planck scale balls (i.e., geodesic balls of radius $\approx \frac{1}{\lambda} $), by a straightforward application of the Continuous Mapping theorem (see \S\ref{sec:band-lim survey}), or, since the asymptotics \eqref{eq:covar band-lim->covar infty} is valid for all $R>0$, tour de force slightly above it, i.e. geodesic balls of radii $R/\lambda$ with $R\rightarrow\infty$ sufficiently slowly.

\vspace{2mm}

In $1$d the arising Gaussian ensemble is that of trigonometric polynomials: It is the Gaussian ensemble of stationary random processes
\begin{equation}
\label{eq:XN rand trig poly}
X_{N}(x)=\sum\limits_{n=1}^{N}(a_{n}\sin(nt)+b_{n}\cos(nt)),
\end{equation}
with $a_{n},b_{n}\sim\Nc(0,1)$ standard i.i.d. Gaussian random variables, $t\in [0,2\pi]$, and $N\ge 1$. The analogous scaling limit process is the Paley-Wiener
process on $\R$, with the covariance $\frac{\sin{x}}{x}$.
In this manuscript we will avoid dealing with the $1$-dimensional case that was resolved ~\cite{GWAJM}, and instead only survey the high dimensional literature.

\subsection{Random spherical harmonics and Arithmetic Random Waves}

The sphere $\Sc^{2}\subseteq\R^{3}$ admits high spectral degeneracy, and, as a result, here, rather than superimposing eigenfunctions belonging to different eigenspaces \eqref{eq:f band lim def},
it allows for the introduction of the ensemble of {\em random spherical harmonics}, that is a particular instance of {\em monochromatic} random waves
(i.e. band-limited functions \eqref{eq:f band lim def} with $\alpha=1$). The Laplace eigenvalues of the sphere are all the numbers $$\lambda_{\ell}^{2}:=\ell(\ell+1),$$
parameterized by $\ell\in\Z_{\ge 0}$ nonnegative integer numbers.
Given $\ell\ge 0$, the eigenspace corresponding to $\lambda_{\ell}^{2}$ is the space of degree-$\ell$ spherical harmonics, of dimension $2\ell+1$, and let
$$\{\eta_{1}=\eta_{\ell,1},\ldots \eta_{2\ell+1}=\eta_{\ell,2\ell+1} \}$$ be its arbitrary $L^{2}$-orthonormal basis. The degree-$\ell$ random spherical harmonic
is
\begin{equation}
\label{eq:Tl spher harm}
T_{\ell}(x) = \frac{1}{\sqrt{2\ell+1}}\sum\limits_{k=1}^{2\ell+1}a_{k}\eta_{k}(x),
\end{equation}
with $a_{k}$ i.i.d. standard Gaussian, and the convenience pre-factor $\frac{1}{\sqrt{2\ell+1}}$ making $T_{\ell}$ univariate. It is easy to generalize the definition of $T_{\ell}(\cdot)$
for higher dimensional spheres. Since the standard multivariate Gaussian distribution is invariant w.r.t. orthogonal transformations, the law of $T_{\ell}$ is invariant w.r.t. the choice of the $L^{2}$-orthonormal basis.

Alternatively (and equivalently) $T_{\ell}(\cdot)$ is the Gaussian random field uniquely prescribed by the covariance function
\begin{equation}
\label{eq:covar Legendre}
r_{\ell}(x,y)=P_{\ell}(\cos {d(x,y)}),
\end{equation}
where $P_{\ell}$ are the
Legendre polynomials parameterized by their degree $\ell\ge 1$, and $d(\cdot,\cdot)$ is the spherical distance between $x,y\in\Sc^{2}$.
For every $\ell\ge 1$, the random field $T_{\ell}(\cdot)$ is isotropic, in the sense that for every rotation $g$ of $\Sc^{2}$, the law of $T(g\cdot)$ is identical to the law
of $T(\cdot)$.

\vspace{2mm}

As all spherical harmonics are either even or odd (w.r.t. to the transformation $x\mapsto -x$ on $\Sc^{2}$), depending on whether $\ell$ is even or odd respectively,
so a.s. is $T_{\ell}(\cdot)$, defined in \eqref{eq:Tl spher harm}. This means that the Gaussian vector $(T_{\ell}(x),T_{\ell}(-x))$ is fully degenerate, hence the Kac-Rice formula
\eqref{eq:KC 2nd mom} below is not applicable on $\Sc^{2}$. On the other hand, it also means that the nodal set of $T_{\ell}$ is invariant w.r.t. $x\mapsto -x$, hence,
to recover all the information on the nodal set of $T_{\ell}$, it is sufficient to restrict
$T_{\ell}$ to a hemisphere, where it is applicable.

\vspace{2mm}

As $\ell\rightarrow\infty$, the asymptotic behaviour of $P_{\ell}$ is described by Hilb's approximation ~\cite[Formula (8.21.17) on p. 197]{Szego}
\begin{equation}
\label{eq:Pl Hilb}
P_{\ell}(\cos\phi) = \left(\frac{\phi}{\sin{\phi}} \right)^{1/2} \cdot J_{0}((\ell+1/2)\phi)+\delta(\phi),
\end{equation}
where $\delta(\cdot)=\delta_{\ell}(\cdot)$ is a {\em small} error term.
As we fix a reference point $x_{0}\in\Sc^{2}$ (e.g. $x_{0}=N$, the northern pole), rescale the variables around $x_{0}$ and flatten the sphere in a small vicinity of $x_{0}$,
the covariance function \eqref{eq:covar Legendre} is approximately
\begin{equation}
\label{eq:covar spher RWM}
r_{\ell}(x,y) \approx  \left(\frac{\psi/\ell}{\sin\left(\psi/\ell\right)} \right)^{1/2} \cdot J_{0}(\psi),
\end{equation}
where $$\psi =\psi(x,y)= (\ell+1/2)\phi= (\ell+1/2)d(x,y)\in [0,(\ell+1/2)\pi],$$ and
$x,y\in\Sc^{2}$, both in the vicinity of $x_{0}$. We identify the factor $J_{0}(\psi)$ in \eqref{eq:covar spher RWM} as the covariance function \eqref{eq:RWM covar} of Berry's RWM;
the extra factor $\left(\frac{\psi/\ell}{\sin\left(\psi/\ell\right)} \right)^{1/2}$, reminiscent of the geometry of the sphere and the rescaled variables, is only close to $1$ if
$d(x,y)\rightarrow 0$, consistent to the general monochromatic waves \eqref{eq:f band lim def} (with $\alpha=1$), though on a bigger neighbourhood, and faster effective convergence rate. Hence the scaled random field $T_{\ell}$, restricted to a small, vanishing, neighbourhood of a reference point (and only restricted one) converges, in appropriate sense,
to the RWM. Mind that the sphere with the geodesic flow is a {\em completely integrable} dynamical system, hence the RWM is not expected to (nor does it)
represent the {\em global} behaviour of the {\em individual} spherical harmonics.

It is known that {\em every} square-summable isotropic random field $G$ on $\Sc^{2}$ can be decomposed in $L^{2}$ as
\begin{equation*}
%G(x)\overset{L^{2}(\Sc^{2})}{=}\sum\limits_{\ell=1}^{\infty}C_{\ell}T_{\ell}(x),
G(x)\overset{L^{2}(\Sc^{2})}{=\joinrel=}\sum\limits_{\ell=1}^{\infty}C_{\ell}T_{\ell}(x)
\end{equation*}
in the sense that
\begin{equation*}
\E\left[\left\|G-\sum\limits_{\ell=1}^{L}C_{\ell}T_{\ell}\right\|_{L^{2}(\Sc^{2})}^{2}\right]\substack{\longrightarrow \\ L\rightarrow\infty} 0.
\end{equation*}
The collection of the nonnegative numbers $\{C_{\ell}\}_{\ell\ge 1}$ is called the {\em power spectrum} of $G$.
That is, $\{T_{\ell}\}_{\ell \ge 1}$ are the Fourier components of every nice Gaussian isotropic field. For example, in cosmology these represent the Cosmic Microwave Background (CMB) radiation measurements under the Gaussianity assumption, and the high energy limit $\ell\rightarrow\infty$ stands for the high precision of those measurements.

Another surface admitting high spectral degeneracy is the standard $2$-dimensional torus $\Tb^{2}=\R^{2}/\Z^{2}$.
Oravecz-Rudnick-Wigman ~\cite{ORW} introduced the ensemble of {\em random toral eigenfunctions}, usually referred to as ``Arithmetic Random Waves" (ARW). Let $$S=\{a^{2}+b^{2}:\:a,b\in \Z\}\subseteq\Z$$ be the set of all positive integers expressible as a sum of two squares, and for $n\in S$ we let $$\Lambda_{n}:= \{\lambda=(\lambda_{1},\lambda_{2})\in\Z^{2}:\: \|\lambda\|^{2}=n\} = \Z^{2}\cap \sqrt{n}\Sc^{1}$$ be
the collection of all lattice points of squared Euclidean norm $n$, or equivalently, lattice points lying on the centred circle of radius $\sqrt{n}$ (alternatively, the nonempty set of representations of $n$ as a sum of two squares). Denote $N_{n}=r_{2}(n):=|\Lambda_{n}|$ to be the {\em number} of lattice points lying on the said circle (also the number of representations of $n$ as sum of two squares). It is well-known that the toral Laplace eigenvalues are all the number of the form $\{4\pi^{2}n:\: n\in S\}$, and the eigenspace corresponding to $4\pi^{2}n$ is the collection of all (complex) linear combinations
of the plane waves $e(\langle \lambda,x\rangle)$, where $x=(x_{1},x_{2})\in\Tb^{2}$, $$ \langle \lambda,\, x\rangle= \lambda_{1}x_{1}+\lambda_{2}x_{2}$$ is the standard Euclidean inner product, and $e(\cdot)=e^{2\pi i \cdot}$.

The Arithmetic Random Waves is the Gaussian ensemble\footnote{Note the abuse of notation, as in this context $n$ denotes the energy rather than the manifold dimension.} $f_{n}:\Tb^{2}\rightarrow\R$ of the random linear combinations of the plane waves
\begin{equation}
\label{eq:fn ARW def}
f_{n}(x) = \frac{1}{\sqrt{2N_{n}}}\sum\limits_{\lambda\in\Lambda_{n}}a_{\lambda}e\left(\langle \lambda,\, x\rangle\right),
\end{equation}
where the coefficients $\{a_{\lambda}\}_{\|\lambda\|^{2}=n}$ are standard complex
Gaussian i.i.d., save for the relation $a_{-\lambda}=\overline{a_{\lambda}}$, making $f_{n}$ real-valued.
As for the spherical harmonics, the convenience pre-factor $\frac{1}{\sqrt{2N_{n}}}$ was introduced to make $f_{n}$ univariate.
Equivalently, the $f_{n}$ are {\em stationary} Gaussian random waves indexed by $\Tb^{2}$, with the covariance function
\begin{equation}
\label{eq:rn ARW covar}
r_{n}(x)=\frac{1}{N_{n}}\sum\limits_{\lambda\in\Lambda_{n}}\cos(2\pi \langle \lambda,x \rangle),
\end{equation}
$x\in \Tb^{2}$. Finally, $f_{n}$ could be uniquely prescribed via (the scaled version of) its spectral measure
supported on $\Sc^{1}\subseteq \R^{2}$:
\begin{equation*}
\mu_{n}:= \frac{1}{N_{n}}\sum\limits_{\lambda\in \Lambda_{n}}\delta_{\lambda/\sqrt{n}},
\end{equation*}
with $\delta_{x}$ the Dirac delta function supported on $x$ (cf. the spectral measure of Berry's RWM described above).

\subsection{Nodal intersections against smooth curves and hypersurfaces}

Let $\M$ be a surface and $\Cc\subseteq\M$ be a smooth curve. One is interested in the following (deterministic) questions: 1.
Is the number of the intersections $M_{\Cc;i}$ of the nodal line of Laplace eigenfunction $\varphi_{i}$ (as in \S\ref{sec:band-lim functions}) against $\Cc$ finite
for $i$ sufficiently large? 2. If the number of nodal intersections is finite, what is its asymptotic behaviour as $i\rightarrow\infty$?
The natural scaling of the problem suggests that, provided that, given $\Cc$, the answer to (1) is ``yes", then the number of nodal intersections
should be asymptotic to $$c_{\Cc}\cdot \sqrt{\lambda_{i}},$$ with an appropriate $c_{\Cc}>0$. Alternatively, at least optimal, up to a constant,
upper and lower bounds $$ c_{\Cc}\cdot \lambda_{i}\le  M_{\Cc;i} \le C_{\Cc}\cdot \lambda_{i},$$ $0<c_{\Cc}<C_{\Cc}$ should be pursued. Of course, same questions
could be posed for manifolds of arbitrarily high dimension, and $\Sigma$, a hypersurface of arbitrary codimension, in place of $\Cc$ (depending on the dimensions of $\M$ and
$\Sigma$, these might not be discrete zeros, whence the relevant quantity is the intersection volume).

These questions were addressed by Toth-Zelditch for complex analytic manifolds, by methods of {\em complexification}, i.e. the number of nodal intersections
is bounded by the number of (complex) zeros of the extension of $\varphi_{i}$ in a complex neighbourhood of $\M$ around $\Cc$. An optimal upper bound was proved ~\cite{ToZe1}
for $\Cc$ satisfying a ``goodness" condition, which, in practice, might be quite difficult to validate for a given $\Cc$, see also the references within ~\cite{ToZe1}. A particularly important scenario is when $\M$ is
a billiard (a planar surface with a boundary), and $\Cc$ is its boundary, addressed in ~\cite{ToZe2}, where an optimal upper bound was proved in a number of scenarios. A precise result for the {\em expected} number of intersections of the band-limited functions against the boundary of a generic billiard was established in ~\cite{ToWig}.

\vspace{2mm}

The same questions for the particular case of (deterministic) {\em toral} eigenfunctions of arbitrary dimension, susceptible to the number theoretic methods, was addressed by Bourgain-Rudnick.
They proved ~\cite{BouRud2,BouRud1} the finiteness and the optimal upper bound for the intersection number (or volume) for generic curves, and a nearly optimal lower bound, along with some other interesting results of relevant nature.

\section{Main ideas: Kac-Rice formulae and their applications}

\subsection{Zero density}

Let $F:\Dc\rightarrow\R$ indexed by the domain $\Dc\subseteq\R^{n}$, be a sufficiently smooth Gaussian random field so that for every $x\in\Dc$, $F(x)$ is non-degenerate,
and, with no loss of generality, one may assume that $F$ is {\em centred}, i.e. for every $x\in\Dc$, $\E[F(x)]=0$.
Assuming that for every $x\in \Dc$, $F(x)$ is nondegenerate Gaussian r.v., one may define the {\em zero density} (also called the {\em first intensity}) of $F$ as
\begin{equation}
\label{eq:K1 zero dens}
K_{1}(x)=K_{F;1}(x):= \phi_{F(x)}(0)\cdot\E[\|\nabla F(x) \|\big| F(x)=0],
\end{equation}
$x\in\Dc$, where $\phi_{F(x)}(0)=\frac{1}{\sqrt{2\pi}\sqrt{\var(F(x))}}$ is the probability density function of $F(x)$ evaluated at $0$, and
$\E[\|\nabla F(x) \|\big| F(x)=0]$ is the Gaussian expectation of the norm of $\nabla F(x)$ conditioned on $F(x)=0$. Then, under very mild assumptions on
$F$, one may evaluate the expected nodal volume of $F$ on compact $\Dc$ as
\begin{equation}
\label{eq:exp zer = int dens}
\E[\vol_{n-1}(F^{-1}(0))] = \int\limits_{\Dc}K_{1}(x)dx.
\end{equation}
To the best knowledge of the author, the mildest known sufficient condition for \eqref{eq:exp zer = int dens} is ~\cite[Theorem 6.8]{AzWs},
applicable under the only extra assumption that $F$ has a.s. no degenerate zeros, for which,
in turn, a generalization of Bulinskaya's Lemma ~\cite[Proposition 6.12]{AzWs} gives a sufficient condition in terms of the boundedness of the density of $(F(x),\nabla F(x))$ at $(0,0)$, uniform w.r.t. $x\in \Dc$.
It is easy to modify the definition of the zero density \eqref{eq:K1 zero dens} corresponding to other local quantities: For example, at level $t\in\R$, it assumes the shape
$$K_{F;1,t}(x):= \phi_{F(x)}(t)\cdot\E[\|\nabla F(x) \|\big| F(x)=t],$$ and it is also an easy (though a bit more complicated) variation to adapt ~\cite{CW,CMW} for the critical points or critical points whose values are restricted to lie in a window. In the Riemannian case the above holds true, except that the meaning of $\nabla F(x)$ is naturally adapted.

If, as we are accustomed, $F$ is stationary, then the density $K_{1}\equiv K_{1}(x)$ does not depend on $x$, and the formula \eqref{eq:exp zer = int dens} assumes a particularly simple shape
\begin{equation*}
\E[\vol_{n-1}(F^{-1}(0))] = \vol(\Dc)\cdot K_{1},
\end{equation*}
i.e. the expected number of the nodal volume of a stationary random field restricted to a domain is proportional to its volume. In this case, the number $K_{1}$ is expressed as a Gaussian integral
\eqref{eq:K1 zero dens}, whose parameters are easily computed in terms of the covariance $r_{F}$ and its mixed derivatives up to $2$nd order, all evaluated at the origin.

\subsection{The $2$-point correlation function}

Further, if for every $x,y\in\Dc$ so that $x\ne y$, the bivariate Gaussian vector $(F(x),F(y))$ is non-degenerate, one can define the {\em $2$-point correlation function} ($2$nd intensity) of the nodal set as
\begin{equation}
\label{eq:K2 def Gauss int}
K_{2}(x,y)=K_{F;2}(x,y):= \phi_{F(x),F(y)}(0,0)\cdot\E[\|\nabla F(x) \|\cdot \|\nabla F(y) \|\big| F(x)=F(y)=0],
\end{equation}
$x,y\in\R^{n}$, $x\ne y$, where $\phi_{F(x),F(y)}(\cdot, \cdot)$ is the Gaussian density of $(F(x),F(y))$, and $$\E[\|\nabla F(x) \|\cdot \|\nabla F(y) \|\big| F(x)=F(y)=0]$$ is the Gaussian
expectation of the product $\|\nabla F(x) \|\cdot \|\nabla F(y) \|$ of norms conditioned on $F(x)=F(y)=0$. For $F$ sufficiently nice, for $n\ge 2$, one may express the $2$nd moment of the nodal volume of $F$ via the double integral
\begin{equation}
\label{eq:KC 2nd mom}
\E[\vol_{n-1}(F^{-1}(0))^{2}] = \int\limits_{\Dc\times\Dc}K_{2}(x,y)dxdy,
\end{equation}
which, for $n=1$ yields the factorial $2$nd moment \footnote{This is the contribution of the so-called ``diagonal" - tuples of the {\em discrete} zeros $(z_{i},z_{i})$ of $F$, see e.g. the argument given by Cramer-Leadbetter ~\cite{CrLe}.}
\begin{equation}
\label{eq:KC 2nd mom fact}
\E[\#(F^{-1}(0))^{2}-\#(F^{-1}(0))] = \int\limits_{[a,b]^{2}}K_{2}(x,y)dxdy
\end{equation}
of the number of the discrete zeros of $F$ on a compact interval $[a,b]$.
To the best knowledge of the author, the mildest sufficient condition for \eqref{eq:KC 2nd mom} (and \eqref{eq:KC 2nd mom fact}) is given by ~\cite[theorems 6.8-6.9]{AzWs}, a little more than the non-degeneracy of $(F(x),F(y))$ for every $x\ne y$. Unless specified otherwise, we will be willing to assume for a while that $n\ge 2$.

The definition \eqref{eq:K2 def Gauss int} of the $2$-point correlation function gives (for $x\ne y$) the value of $K_{2}(x,y)$ in terms of a certain {\em Gaussian integral} depending on the covariance function $r_{F}$ and its mixed derivatives up to $2$nd order, evaluated at $(x,y)$ and the diagonal points $(x,x)$ and $(y,y)$; in $1$d one may evaluate that integral explicitly to express it as an elementary function of the said derivatives.
If $F$ is stationary, then one may normalize it to be univariate: for every $x\in\Dc$, $\var(F(x))=1$. Then (assuming $n\ge 2$),
the $2$-point correlation function depends on $z:=y-x$ only, and it simplifies to\footnote{We tacitly assume that $0\in\Dc$, otherwise shift the random field.}
\begin{equation*}
K_{2}(z)=\frac{1}{2\pi\sqrt{1-r_{F}(z)^{2}}}\E[\|\nabla F(0) \|\cdot \|\nabla F(z) \|\big| F(0)=F(z)=0],
\end{equation*}
whereas \eqref{eq:KC 2nd mom} reads
\begin{equation}
\label{eq:KC 2nd mom stat}
\E[\vol_{n-1}(F^{-1}(0))^{2}] = \int\limits_{\Dc\times\Dc}K_{2}(x-y)dxdy,
\end{equation}
with an obvious adaptation for $n=1$. When $F$ is defined on a manifold possessing some symmetries w.r.t. which the law of $F$ is invariant (such as the standard torus with respect to
translations, or the sphere with respect to rotations), the double integral on the r.h.s. of \eqref{eq:KC 2nd mom stat} is expressible in terms of a simple integral.

\vspace{2mm}

Typically, for $F$ stationary, defined on the whole of $\R^{n}$, the asymptotic behaviour of $K_{2}(z)$ as $|z|\rightarrow\infty$ determines the asymptotic law for the variance of the nodal volume
$\vol_{n-1}(F^{-1}(0)\cap B(R))$ of $F$ restricted to increasing balls $B(R)$, $R\rightarrow\infty$; in the Riemannian setting with scaling, the relevant asymptotic behaviour is for the {\em scaled} variables diverging.
However, the $2$-point correlation function carries so much more information than merely a precise evaluation of the $2$nd moment $\E[\vol_{n-1}(F^{-1}(0)\cap \Dc)^{2}]$ or the variance $\var(\vol_{n-1}(F^{-1}(0)))$, or the asymptotic nodal volume variance restricted to increasing balls. First, it gives the variance on arbitrary shapes or their homotheties with no re-calculation, and in the Riemannian setting, the variance of the nodal volume of $f_{k}$ restricted to subdomains, possibly shrinking slower than the scaling (see the applications given in \S\ref{sec:nod length} and \S\ref{sec:gen nonlin func}). It also endows the notions of zero attraction or repulsion with a proper meaning: we say that if zeros repel if $K_{2}(x,y)$ vanishes as $\|x-y\|\rightarrow 0$, and attract if $K_{2}(x,y)$ grows to infinity as $\|x-y\|\rightarrow 0$.

When in addition, $F$ is {\em isotropic}, then $K_{2}(x,y)= K_{2}(\|x-y\|)$ depends on the distance between $x$ and $y$ only. Hence, in this case, one encodes the attraction of repulsion of zeros in terms of the value of
$K_{2}(0)$ at the origin, understood in the limit sense; for the former (attraction) $K_{2}(0)=+\infty$, whereas for the latter, $K(0)=0$. In the isotropic case, one may also entertain an intermediate notion of zero attraction or repulsion, in the situation when $K_{2}(0)>0$ is a finite, strictly positive number. For example, one may compare ~\cite{BCW1,BCW2} this value at the origin of the $2$-point correlation function corresponding to the critical points of $F$ to that of the Poisson point process of the same intensity.

\vspace{2mm}

Analogously to the above, one may introduce the $k$-point correlation function, $k\ge 3$, to relate to the higher moments of the nodal volume.
It is easy to modify the definition of the zero density \eqref{eq:K1 zero dens} and the other correlation functions for other local quantities: For example, at level $t\in\R$, the density function assumes the shape
$$K_{F;1,t}(x):= \phi_{F(x)}(t)\cdot\E[\|\nabla F(x) \|\big| F(x)=t],$$ and it is also an easy (though a bit more complicated) variation to adapt ~\cite{CW,CMW} for the critical points or critical points whose values are restricted to lie in a window. In the Riemannian case all of the above holds true, except that the meaning of $\nabla F(x)$ is naturally adapted.

\subsection{Approximate and Mixed Kac-Rice formulae}
\label{sec:Approx Mixed KR}

Provided that the said sufficient conditions on $F$ are satisfied, the Kac-Rice integral \eqref{eq:KC 2nd mom} on $\Dc$ computes the second moment (or the second factorial moment in $1$d) on the whole of $\Dc$, and hence its variance\footnote{For $n=1$ obvious adjustments are due.}
\begin{equation}
\label{eq:var=K2-K1K1}
\begin{split}
\var\left(\vol_{n-1}(F^{-1}(0))\right) &=\int\limits_{\Dc\times\Dc}K_{2}(x,y)dxdy-(\E[\vol_{n-1}(F^{-1}(0))])^{2}\\&= \int\limits_{\Dc\times\Dc}\left(K_{2}(x,y)-K_{1}(x)K_{1}(y)\right)dxdy,
\end{split}
\end{equation}
$n\ge 2$. If $\Dc_{1},\Dc_{2}\subseteq \Dc$ are two domains, then one may obtain the {\em covariance} of the nodal volume of $F$ restricted to $\Dc_{1}$ and $\Dc_{2}$ respectively by
restricting the range of the integral on the r.h.s. of \eqref{eq:var=K2-K1K1}:
\begin{equation}
\label{eq:cov=K2-K1K1}
\begin{split}
&\cov\left(\vol_{n-1}(F^{-1}(0)\cap \Dc_{1}),\vol_{n-1}(F^{-1}(0)\cap \Dc_{2}) \right) \\&=
\int\limits_{\Dc_{1}\times\Dc_{2}}K_{2}(x,y)dxdy-\E[\vol_{n-1}(F^{-1}(0)\cap \Dc_{1})]\cdot \E[\vol_{n-1}(F^{-1}(0)\cap \Dc_{2})]\\&= \int\limits_{\Dc_{1}\times\Dc_{2}}\left(K_{2}(x,y)-K_{1}(x)K_{1}(y)\right)dxdy.
\end{split}
\end{equation}
The upshot is that, in some important cases, when \eqref{eq:var=K2-K1K1} either does not hold, or it is technically demanding to validate its sufficient conditions, \eqref{eq:cov=K2-K1K1} may still
be verifiable for a wide family of carefully chosen pairs of non-intersecting subdomains $\Dc_{1},\Dc_{2}\subseteq\Dc$. This way it might still be possible to justify that \eqref{eq:var=K2-K1K1} (``Approximate Kac-Rice" ~\cite[Proposition 1.3]{RuWiAJM}, see \S\ref{sec:ARW ni} below), or, in some cases, further approximation of \eqref{eq:var=K2-K1K1}, holds {\em asymptotically} ~\cite{RuWiYe} in some regime,
also allowing for the asymptotic analysis of the variance $\var\left(\vol_{n-1}(F^{-1}(0))\right)$ of the {\em total} nodal volume.

\vspace{2mm}

Let $F,G:\Dc\rightarrow\R$, $\Dc\subseteq\R^{n}$ with\footnote{For $n=1$ some adjustment is required. Otherwise, it is applicable away from the diagonal unimpaired.} $n\ge 2$ (or $F,G$ defined on $\M$)
be {\em two} Gaussian random fields defined on the same probability space. In analogy to the above, under suitable conditions on $F,G$, one may define the {\em cross-correlation} function (see e.g.
~\cite[page $37$]{GWAJM} or ~\cite[p. 7-8]{MRW}; it is contained in the $1970$'s literature in the $1$d case)
\begin{equation*}
\widetilde{K_{2}}(x,y)=\widetilde{K}_{2;F,G}(x,y):= \phi_{F(x),G(y)}(0,0)\cdot\E[\|\nabla F(x) \|\cdot \|\nabla F(y) \|\big| F(x)=F(y)=0],
\end{equation*}
expressible in terms of the cross-covariance $$\widetilde{r}_{F,G}(x,y) = \E[F(x)\cdot G(y)]$$ and its mixed derivatives up to $2$nd order.
Under suitable conditions, the corresponding (``mixed") Kac-Rice integral computes the {\em covariance} between the zeros of $F$ and $G$ on\footnote{Similarly to above, one may restrict the range of the integral to $\Dc_{1}\times\Dc_{2}$, $\Dc_{1},\Dc_{2}\subseteq\Dc$, resulting in the covariance between the nodal volume of $F$ and $G$ restricted to $\Dc_{1}$ and $\Dc_{2}$ respectively.} $\Dc$:
\begin{equation}
\label{eq:cov mixed KR}
\cov(\vol_{n-1}(F^{-1}(0)),\vol_{n-1}(G^{-1}(0))) = \int\limits_{\Dc\times\Dc}\left(\widetilde{K}_{2}(x,y) - K_{1;F}(x)\cdot K_{1;G}(y)\right)dxdy.
\end{equation}

The Mixed Kac-Rice formula is in particular useful when either $F$ represents an ensemble $\{f_{k}\}$ of Gaussian random fields converging in law to a limit Gaussian random field $G$, or when $G$ is a mollified version of $F$; for example, $G$ could be stationary or even isotropic, whereas $F$ is only asymptotically such, or $G$ could be $M$-dependent\footnote{That is, $r_{G}(x,y)=0$ for $\|x-y\|>M$.} for $M$ sufficiently slowly growing parameter. This way, from quantitative convergence of $r_{F}(\cdot,\cdot)$ and its mixed derivatives to $r_{G}(\cdot,\cdot)$ and its derivatives respectively, one may infer various properties on $F$ from the analogue properties on $G$: for example, one may infer ~\cite{GWAJM} the variance and the limit law of the total number of zeros of random trigonometric polynomials \eqref{eq:XN rand trig poly} from the Paley-Wiener process on $\R$. The proximity alone of $r_{F}$ to $r_{G}$ does not guarantee the success of this procedure; the key step is constructing a coupling of $F$ and $G$ so that the cross-covariance $\widetilde{r}_{F,G}$ is quantitatively close, in the suitable norm, to $r_{F}$ (and $r_{G}$). Combining the two ideas of Approximate and Mixed Kac-Rice to express covariances of nodal volume on different domains can give extra boost to these powerful methods.
As another example, for $n\ge 2$, one may reproduce this argument to infer the variance and the limit law of the nodal volume of band-limited functions \eqref{eq:f band lim def} from the corresponding limit random field, in somewhat more restrictive scenario, only in geodesic balls of radius\footnote{This is a by-product of the fact that the convergence of $r_{f_{\lambda}}$ to the limit covariance is restricted only slightly above Plack scale
if $\M$ has no conjugate points; otherwise, other than Zoll surfaces or rational tori, no general result above Plack scale is known.} logarithm power higher than the Plack scale.

\vspace{2mm}

Rather than comparing the nodal volume of two different random fields, as in \eqref{eq:cov mixed KR}, one may derive the cross-correlation formula that compares between {\em different} properties of the same random field, or
different random fields, and integrate it to obtain the covariance between these. For example, this way one may compute the covariance between the nodal volume and the number of critical points of the same random field $F$
on a domain $\Dc$ or, for instance, the covariance of the volumes of two nonzero levels $\cov(\vol_{n-1}(F^{-1}(t_{1})),\vol_{n-1}(F^{-1}(t_{2})))$, or, more generally,
the covariance of any local properties of two random fields restricted to arbitrary domains.

\subsection{Auxiliary function for linear statistics, or Euclidean random fields}
\label{sec:W transform}

In a situation when $F:\R^{2}\rightarrow \R$ is a stationary {\em isotropic} Gaussian process, for example, Berry's RWM (alternatively, $F:\Sc^{2}\rightarrow\R$ is invariant w.r.t. rotations), the $2$-point correlation function $$K_{F;2}(x,y)=K_{2}(\|x-y\|)$$ is a function of the distance $\|x-y\|$ (resp. of the spherical distance $d(x,y)$). One may then use \eqref{eq:KC 2nd mom} with Fubini to rewrite
\begin{equation*}
\E[\vol_{n-1}(F^{-1}(0))^{2}] = \int\limits_{0}^{\diam(\Dc)}K_{2}(t)W_{\Dc}(t)dt,
\end{equation*}
where $\diam(\Dc)$ is the diameter of $\Dc$, and $W_{\Dc}(\cdot)$ is the auxiliary function ~\cite{W fluctuations}
\begin{equation}
\label{eq:W tranform def bas}
W_{\Dc}(t)=\meas(\{(x,y)\in \Dc^{2}:\: \|x-y\|=t \}) = \int\limits_{\Dc}\len(\{y\in\Dc:\: \|x-y\|=t\})dx,
\end{equation}
independent of $F$.
Given $\Dc$ finding the precise values of $W_{\Dc}(\cdot)$ might be quite technically demanding or not have a closed formula, even in some elementary cases, e.g. $\Dc$ is a disc. However, in many cases, such as the situation of the random spherical harmonics, the bulk of the variance \eqref{eq:KC 2nd mom} is concentrated around the diagonal, hence what matters is only the asymptotics of $W_{\Dc}(\cdot)$ near the origin.

More generally, one may consider the {\em linear statistics}
\begin{equation*}
\vol^{\varphi}(F):= \int\limits_{F^{-1}(0)}\varphi(x)dx,
\end{equation*}
where $\varphi:\R^{2}\rightarrow\R$ (or $\varphi:\M\rightarrow\R$) is some test function, e.g. it could be a smooth function or the characteristic function $\varphi=\mathds{1}_{\Dc}$ of a domain $\Dc\subseteq\R^{2}$
or $\Dc\subseteq\M$. Then one has ~\cite{W fluctuations}\footnote{It is easy to
combine it with the ideas presented towards the end of \S\ref{sec:Approx Mixed KR},
to express $$\cov(\vol^{\varphi_{1}}(F),\vol^{\varphi_{2}}(G))$$ for test functions $\varphi_{1},\varphi_{2}$, and stationary isotropic Gaussian random fields $F,G$.
It is also easy to generalize it for non-isotropic random fields.}
\begin{equation*}
\E[(\vol^{\varphi}(F))^{2}] = \int\limits_{0}^{\diam(\Dc)}K_{2}(t)W_{\varphi}(t)dt
\end{equation*}
with
\begin{equation*}
W_{\varphi}(t)=\iint\limits_{(x,y)\in\R^{2}\times\R^{2}:\: \|x-y\|=t}\varphi(x)\varphi(y)dxdy = \int\limits_{\R^{2}}dx\int\limits_{\partial B_{x}(t)}\varphi(y)dy.
\end{equation*}
The function $W_{\varphi}$ has a {\em smoothing} effect on $\varphi$ in a way similar to convolving functions; in fact, it is possible to express it as a convolution of certain by-products of $\varphi$.
It gives a powerful tool for evaluating the variance of linear statistics (in particular, of the nodal volume restricted to a subdomain of an ambient manifold) of the nodal lines, assuming that the asymptotic
behaviour of $K_{2}(\cdot)$ is understood, with no need for re-computations or solving the difficult geometric problem \eqref{eq:W tranform def bas}, that is independent of $F$.

\subsection{Wiener Chaos expansion of Kac-Rice integrals}

For a ``nice" Gaussian random field $F:\R^{n}\rightarrow\R$ (alternatively, $F:B(R)\rightarrow\R$ or $F:\M\rightarrow\R$ etc.), defined on the sample space $\Omega$, let
$\mathcal{X}=\mathcal{X}_{F}\subseteq L^{2}(\Omega)$ be the (complex) closed Gaussian Hilbert space generated by the vectors $\{F(x):\: x\in\R^{n}\}$, and $\Hc=\Hc_{F}\subseteq L^{2}(\Omega)$ be the Hilbert space of random variables in $L^{2}(\Omega)$, measurable w.r.t. the sub-$\sigma$-algebra generated by $\mathcal{X}$. For $q\ge 0$ the $q$'th Wiener chaos space is the closed span $\Hc_{q}$ of
$$\{H_{\alpha}(\xi_{1},\ldots,\xi_{k}):\: k\ge 1, \, \alpha=(\alpha_{1},\ldots \alpha_{k})\in\Z_{\ge}^{k},\,\alpha_{1}+\ldots\alpha_{k}=q  \},$$ where the $H_{\alpha}(\cdot)$ is the $k$-variate Hermite polynomial of multi-index $\alpha$.
For example, $\Hc_{0}$ consists of all a.s. constant random variables.
One has the Wiener chaos decomposition
\begin{equation}
\label{eq:Wiener chaos decomp}
\Hc = \bigoplus\limits_{q=0}^{\infty}\Hc_{q} ,
\end{equation}
orthogonal w.r.t. the inner product $$\langle\xi,\xi'\rangle=\E[\xi\cdot\overline{\xi}],$$ that is, if $\xi\in\Hc_{q}$, $\xi'\in \Hc_{q'}$ with $q\ne q'$, then $\xi,\,\eta$ are {\em uncorrelated}.

Using Federer's Co-Area formula ~\cite{Federer}, one may formally express the nodal volume of a sufficiently smooth function $f:\M\rightarrow\R$ as
\begin{equation}
\label{eq:coarea delta}
\vol_{n-1}(f^{-1}(0)) = \int\limits_{\M}\delta_{0}(f(x))\cdot \|\nabla f(x)\|dx,
\end{equation}
where $\delta$ is the Dirac delta function\footnote{In reality, one replaces $\delta_{0}(\cdot)$ with $\frac{1}{2\epsilon}\chi_{[-\epsilon,\epsilon]}(\cdot)$, performs all the steps
below, and takes $\epsilon \downarrow 0$.}. More generally, if $g:\M\rightarrow\R^{k}$, where $\M$ is a smooth $n$-manifold, and $k\le n$, then the nodal $(n-k)$-volume of $g$ is given by
\begin{equation*}
\vol_{n-k}(g^{-1}(0)) = \int\limits_{\M}\delta_{0}(f(x))\cdot \|J_{g}(x)\|dx,
\end{equation*}
where $J_{g}(\cdot)$ is the volume of the parallelogram spanned by the row vectors of the Jacobi matrix of $g$, e.g., if $n=k$, then $J_{g}(\cdot)$ is the Jacobian of $g$. One can infer from
the expression \eqref{eq:coarea delta}, that, given
a random field $F:\M\rightarrow\R$ (or, for example $F:B(R)\rightarrow\R$ with the Euclidean ball $B(R)\subseteq\R^{n}$)
the random variable $\Vc:=\vol_{n-1}(f^{-1}(0))$ is measurable w.r.t. the random variables generated by $\mathcal{X}$, hence it satisfies $\Vc\in \Hc$, if only it has finite variance,
whence one may project it $$\Vc_{q}:=\proj_{\Hc_{q}}(\Vc)$$
onto $\Hc_{q}$ for $q\ge 0$, and thus decompose it w.r.t. the decomposition \eqref{eq:Wiener chaos decomp}.
For example, the $0$'th Wiener chaotic component of $\Vc$
\begin{equation*}
\Vc_{0}=\E[\Vc]
\end{equation*}
is its expectation, as one has for an arbitrary $X\in \Hc$.

This way, rather than study the distribution of the stochastic integral of the type \eqref{eq:coarea delta}, one may decompose it into the components $\Vc_{q}$, as above, and study their distributions separately,
before combining the individual results for the purpose of obtaining a limit law of the nodal volume $\Vc$. This was first done by Kratz-Le\'{o}n ~\cite{KrLe} for Euclidean random fields $F:\R^{n}\rightarrow\R$ over expanding regions, such as $B(R)$, $R\rightarrow\infty$. This approach was used in ~\cite{MaWi3,MaWi1,MaWi2} (and then elsewhere) in order to study the generic functionals of the random spherical harmonics \eqref{eq:Tl spher harm}, using a somewhat more sophisticated modern language developed in the mean time, that saves one from approximating the given random fields by $M$-dependent ones (see \S\ref{sec:gen nonlin func} below for some details).

\section{Survey of results: Band-limited functions, random spherical harmonics and their Euclidean limits}

\subsection{Berry's random wave model}

Recall that Berry's random wave model is the centred isotropic Gaussian random field $u:\R^{2}\rightarrow\R$ whose covariance is $$r_{u}(x,y)=r_{u}(x-y)=J_{0}(\| x-y\|).$$
By the standard asymptotics of the Bessel functions, as $\psi\rightarrow\infty$, $J_{0}(\psi)$ admits the $2$-term asymptotics\footnote{A full asymptotic expansion is well-known, $2$ terms are sufficient
for the presented results.}
\begin{equation}
\label{eq:J0 2 term asymp}
J_{0}(\psi) = \sqrt{\frac{2}{\pi}} \frac{\sin(\psi+\pi/4)}{\psi^{1/2}}- \frac{1}{8}\sqrt{\frac{2}{\pi}} \frac{\cos(\psi+\pi/4)}{\psi^{3/2}}+O\left(  \frac{1}{\psi^{5/2}} \right),
\end{equation}
yielding the asymptotics as $\|x-y\|\rightarrow\infty$ for $r_{u}(x,y)$, and one may also derive. Alternatively, one may derive the asymptotic expansion of $J_{0}(\psi)$ by writing its expression as the Fourier
transform of the characteristic function of the unit circle, and using the {\em stationary phase} method.

Since $u$ is stationary (ditto isotropic), its zero density is constant on $\R^{2}$, easily evaluated to be\footnote{Though ~\cite{Berry02} did not state it in this precise form, it is easy to extract it from Berry's
manuscript ~\cite{Berry02}.} ~\cite{Berry02}
\begin{equation}
\label{eq:K1 RWM}
K_{1}(x)=K_{u;1}(x)\equiv \frac{1}{2^{3/2}},
\end{equation}
hence the expected nodal length of $u$ restricted to $B(R)$ equals {\em precisely} to
\begin{equation*}
\E\left[\len(u^{-1}(0)\cap B(R))\right] = \frac{1}{2^{3/2}}\cdot \pi R^{2}.
\end{equation*}
Berry then exploited the asymptotics \eqref{eq:J0 2 term asymp} of the covariance function and its mixed derivative to yield a
$2$-term asymptotic expansion for the $2$-point correlation function at infinity: Since $u$ is isotropic, $K_{u;2}(x,y)=K_{2}(\|x-y\|)=K_{2}(\psi)$ depends only on the distance $\psi:=\|x-y\|$
between the two points, and he asserted that
\begin{equation}
\label{eq:Berry's 2 term cancellation}
K_{2}(\psi) = \frac{1}{8}+\frac{\sin(2\psi)}{4\pi \psi}+\frac{1}{512\pi^{2}}\frac{1}{\psi^{2}} +E(\psi),
\end{equation}
where $E(\psi)$ contains {\em purely oscillatory}
terms of order of magnitude $\frac{1}{\psi^{2}}$ as well as a term of order of magnitude $O\left( \frac{1}{\psi^{3}}  \right)$.

Berry's asymptotics \eqref{eq:Berry's 2 term cancellation} is quite remarkable, in more than one way. First, it gives a precise asymptotic law for the variance of the nodal volume of $u$ restricted to $B(R)$
as $R\rightarrow\infty$: only the regime $\|x-y\|\rightarrow\infty$ contributes to the integral on the r.h.s. of \eqref{eq:var=K2-K1K1}, the constant term $\frac{1}{8}$ in
\eqref{eq:Berry's 2 term cancellation} perfectly cancels out with the squared expectation via \eqref{eq:K1 RWM}, since at infinity (i.e. $\|x-y\|\rightarrow\infty$),
$(u(x),\nabla u(x))$ and $(u(y),\nabla u(y))$ are stochastically independent. Further, one needs to integrate the two remaining terms on the r.h.s. of \eqref{eq:Berry's 2 term cancellation} on $B(R)\times B(R)$
(recall that $\psi:=\|x-y\|$), on dismissing the error term $E(\cdot)$. Berry argued that, given a ``generic" $x\in B(R)$ and $\psi \in (0,R)$, for the purpose of deriving the leading asymptotics of the double integral,
one may assume that $y$ is unrestricted, lying on the circle of radius $\psi$ centred at $x$. This simplified model makes the integral into simple, $1$-dimensional integral, for it disregards the geometric factor that is naturally induced on us when we are to evaluate the measure of the set of tuples $$\{(x,y)\in B(R)\times B(R):\: \|x-y\|=\psi\}$$ at distance $\psi$, that, to the author's best knowledge, does not have a simple or elegant answer (cf. \S\ref{sec:W transform}).

The other remarkable fact about Berry's asymptotics \eqref{eq:Berry's 2 term cancellation} is that the remaining leading term $\frac{\sin(2\pi)}{4\pi \psi}$ is purely oscillatory, and will not contribute to the integral \eqref{eq:var=K2-K1K1}. Instead, the bulk of the variance comes from the term $\frac{1}{512\pi^{2}}\frac{1}{\psi^{2}}$, which gives a {\em logarithmic} contribution, and hence the variance is
\begin{equation}
\label{eq:Berry log}
\var(u^{-1}(0)\cap B(R))\sim\frac{1}{128}R^{2}\log{R},
\end{equation}
of lower order than $\sim *R^{3}$, that would be expected by the natural scaling considerations of the problem.
This remarkable phenomenon was named ~\cite{W fluctuations} ``Berry's cancellation", that is highly acclaimed in the recent literature. It has many different, easily identifiable, appearances in other contexts,
see e.g. \S\ref{sec:res spher harm} or \S\ref{sec:survey ARW} below. Though in his original manuscript ~\cite{Berry02}, Berry did not pursue a complete rigorous proof for the asymptotics \eqref{eq:Berry log} for the nodal length variance, it could be fully validated by adding some details: Bound for the contribution of the origin and the intermediate regimes, more detailed derivation of \eqref{eq:Berry's 2 term cancellation}, and the ``symmetrization" of the domain; all of this was performed for the random spherical harmonics (cf. \S\ref{sec:res spher harm}), with the understanding that it could have been fully performed by Berry had he pursued a rigorous proof for \eqref{eq:Berry log}.
Interestingly, Berry's cancellation does not persists ~\cite{W Dartmouth} for nonzero levels, see an explanation in \S\ref{sec:res spher harm}\footnote{In ~\cite{W Dartmouth} it concerned the same situations for random spherical harmonics \eqref{eq:Tl spher harm}. But, as explained in \S\ref{sec:res spher harm} below, given the result for random spherical harmonics, the corresponding result for Berry's RWM is automatic.}.

Berry's original elegant derivation of the asymptotics \eqref{eq:Berry's 2 term cancellation} capitalized on the particular structure of the Gaussian integral \eqref{eq:K2 def Gauss int}.
Since the covariance function $r_{u}$
(and all of its derivatives) decay at infinity, the covariance matrix of $(u(x),u(y),\nabla u(x),\nabla u(y))$ (since $u$ is isotropic, $K_{2}$ is invariant w.r.t. translations and rotations in $\R^{2}$, so, accordingly,
we may assume that $x=(0,0), y=(\psi,0)$), properly normalized, is a perturbation of the identity matrix, hence the resulting Gaussian integral approximates one with separation of variables, whose value cancels the squared expectation. In general (for example, when treating the Gaussian integral corresponding to counting critical points rather than the nodal volume ~\cite{CW}), one invokes
perturbation theory\footnote{It is also the Taylor expansion of the Gaussian integral as a function of the (perturbed) covariance matrix, about the (unperturbed) identity matrix, that we get at infinity.} (i.e. study the eigenvalues and eigenfunctions of the perturbed covariance matrix) to derive an asymptotic formula for $K_{2}(x,y)$ as $|x-y|\rightarrow\infty$. This approach also shows that,
if one is interested in an upper bound for the variance only, then one may control it via the $L^{2}$-norm of the covariance function and its derivatives, properly normalized, for example, via a bound of the type
\begin{equation}
\label{eq:var bnd L2 enrgy covar}
\begin{split}
\var(\len (u^{-1}(0)\cap B(R))) &\ll \|r_{u}\|_{L^{2}(B(R)\times B(R))}^{2}+\|\partial_{x}r_{u}\|_{L^{2}(B(R)\times B(R))}^{2}+
\|\partial_{y}r_{u}\|_{L^{2}(B(R)\times B(R))}^{2}\\&+\|\partial_{xx}r_{u}\|_{L^{2}(B(R)\times B(R))}^{2}+\|\partial_{xy}r_{u}\|_{L^{2}(B(R)\times B(R))}^{2}+\|\partial_{yy}r_{u}\|_{L^{2}(B(R)\times B(R))}^{2}.
\end{split}
\end{equation}

\subsection{Random spherical harmonics: nodal volume and generic functionals}
\label{sec:res spher harm}

%%%%%%%%%%%%%%%%%%%%%%%%%%%%%%%%%%%%Anything holding for RSH holds for RWM, but not vice versa

\subsubsection{Nodal volume}
\label{sec:nodal volume}

Let us now consider the analogue problem for the random spherical harmonics \eqref{eq:Tl spher harm}. Here both the distribution of the {\em total} nodal length, and the nodal length of the function restricted to subdomains, possibly shrinking as $\ell \rightarrow\infty$, makes sense. The asymptotics \eqref{eq:covar spher RWM}
shows that $T_{\ell}$, restricted to spherical caps of shrinking radii (e.g. radius $1/\log{\ell}$), and rescaled by a factor of $\ell$ (rather, $\ell+1/2$), defined on the flattened coordinates of the sphere, converges to Berry's RWM. Since the law of $T_{\ell}(\cdot)$ is invariant w.r.t. rotations of the sphere, it follows that the zero density is independent of $x\in\Sc^{2}$, hence only depends on $\ell$. An easy and straightforward computations shows that
\begin{equation}
\label{eq:K1 spher harm}
K_{\ell;1}(x)\equiv \frac{1}{2^{3/2}}\sqrt{\lambda_{\ell}},
\end{equation}
in particular, the expected nodal length of $T_{\ell}$ is given precisely by
\begin{equation}
\label{eq:exp nod length Tl}
\E[\len(T_{\ell}^{-1}(0))] = \sqrt{2}\pi\cdot \sqrt{\lambda_{\ell}},
\end{equation}
established by B\'{e}rard ~\cite{Berard} by different methods.
Neuheisel ~\cite{Neuheisel} gave the useful upper bound
\begin{equation*}
\var(\len(T_{\ell}^{-1}(0))) = O\left( \ell^{2-1/7} \right),
\end{equation*}
showing that, with high probability, the nodal length is concentrated around the mean. Neuhesel's bound was improved ~\cite{WJMP} to
\begin{equation*}
\var(\len(T_{\ell}^{-1}(0))) = O\left( \ell^{2-1/2} \right).
\end{equation*}

The question of asymptotic variance was finally resolved in ~\cite{W fluctuations}. It was shown that, up to an admissible error term depending on $\ell$, an analogue of \eqref{eq:Berry's 2 term cancellation}
holds in this case, namely, that\footnote{The normalization in ~\cite{W fluctuations} was slightly different, and here we present the result so that to compare it to \eqref{eq:Berry's 2 term cancellation}.}
for $x,y$ belonging to a (fixed) hemisphere
\begin{equation}
\label{eq:K2 sphere expand}
K_{2;\ell}(x,y) = \ell\cdot \left(  \frac{1}{8}+\frac{1}{4\pi\ell}\cdot \frac{\sin(2\psi)}{\sin(\psi/(\ell+1/2))} + \frac{1}{512\pi^{2}\ell}\cdot  \frac{\sin(2\psi)}{\psi\cdot \sin(\psi/(\ell+1/2))}    +E_{\ell}(\psi)  \right),
\end{equation}
with, as above, $\psi=(\ell+1/2)d(x,y)$. Integrating \eqref{eq:K2 sphere expand} on a hemisphere yields the precise asymptotics ~\cite[Theorem 1.1]{W fluctuations}
\begin{equation}
\label{eq:nodal length variance Tl}
\var(\len(T_{\ell}^{-1}(0))) = \frac{1}{32}\log{\ell}+O(1)
\end{equation}
for the nodal length variance of $T_{\ell}$.

A few observations are due. In accordance with the locality of the convergence of the scaled and flattened $T_{\ell}(\cdot)$ to Berry's RWM, so is the convergence of \eqref{eq:K2 sphere expand} to
\eqref{eq:Berry's 2 term cancellation}, e.g. not holding for $\psi\approx \ell$. However, since the bulk of the variance is contained within a small neighbourhood of the diagonal, e.g. inside
$$\{(x,y):\: d(x,y)<1/\log{\ell}\},$$ where, in addition, the sphere is close to be flat, the asymptotic law of the variance \eqref{eq:nodal length variance Tl} is the same as \eqref{eq:Berry log}, properly scaled and understood. One may however notice that, if one scales the sphere by $\ell$, then there is a discrepancy of factor $2$ between the leading constants of \eqref{eq:nodal length variance Tl} compared to \eqref{eq:Berry log} per unit area. This is a by-product of the fact that the nodal line of $T_{\ell}$ is fully determined by its restriction to a hemisphere, hence the nodal length is twice the nodal length of $T_{\ell}$ restricted to a hemisphere.

\vspace{2mm}

In the same paper ~\cite{W fluctuations}, the result \eqref{eq:nodal length variance Tl} was extended in two directions. Let $\varphi:\Sc^{2}\rightarrow\R$ be a $C^{1}$-{\em smooth} test function, and consider the {\em linear statistics}
\begin{equation*}
\Vc_{\ell}^{\varphi}=\Vc^{\varphi}( T_{\ell}):=\int\limits_{T_{\ell}^{-1}(0)} \varphi(x)dx
\end{equation*}
of the nodal line of $T_{\ell}$. For example, if $\varphi=\mathds{1}_{A}$ is the characteristic function of a nice subdomain $A\subseteq \Sc^{2}$, then
$\Vc_{\ell}^{\varphi}$ is identified as simply the nodal length of $T_{\ell}$ restricted to $A$ (though, in this case, $\varphi$ is not smooth). To avoid tautologies, one may assume w.l.o.g. that $\varphi$ is
an {\em even}, not identically vanishing, function. Using \eqref{eq:K1 spher harm} it is easy to evaluate the expectation
\begin{equation*}
\E[\Vc_{\ell}^{\varphi}] = c_{0}(\varphi)\cdot \sqrt{\lambda_{\ell}}.
\end{equation*}
One may further exploit \eqref{eq:K2 sphere expand} with the auxiliary function $W^{\varphi}$ is in \S\ref{sec:W transform}
to obtain the asymptotic law for the variance ~\cite[Theorem 1.4]{W fluctuations}
\begin{equation}
\label{eq:var lin stat}
\E[\Vc_{\ell}^{\varphi}] = c_{1}(\varphi)\cdot \log{\ell}+O_{\varphi}(1),
\end{equation}
with some $c_{1}(\varphi)>0$. Only the behaviour of $W^{\varphi}$ at the origin will determine the variance \eqref{eq:var lin stat}.

Since the characteristic function of a subdomain of $\Sc^{2}$ is not continuous, let alone $C^{1}$, the result \eqref{eq:var lin stat} does not apply in this case. Instead,
the result \eqref{eq:var lin stat} was also proved ~\cite[Theorem 1.5]{W fluctuations} for a different class of functions, namely for bounded functions of {\em bounded variation} that do contain the characteristic functions
of nice subdomains of the sphere. To prove it, a density argument is given, with the use of \eqref{eq:var lin stat} for smooth $\varphi$ approximating the given characteristic function,
and effective control over the dependence of the error term on $\varphi$ as the basis.

\vspace{2mm}

Importantly, only the nodal case, invariant w.r.t. multiplication of $T_{\ell}$ by any real number,
is susceptible to Berry's cancellation, whereas the nonzero levels are not. Let $t\in\R\setminus \{0\}$, and consider the level curve length
$\len(T_{\ell}^{-1}(t))$. Then it is easy to evaluate the density in this case to yield the expectation
\begin{equation*}
\E\left[\len(T_{\ell}^{-1}(t))\right] = e^{-t^{2}/2}\cdot \sqrt{2}\pi\cdot \sqrt{\lambda_{\ell}},
\end{equation*}
cf. \eqref{eq:exp nod length Tl}. For the variance one has ~\cite[Formula (18)]{W Dartmouth}
\begin{equation}
\label{eq:var level set}
\var\left(\len(T_{\ell}^{-1}(t))\right) = c\cdot t^{4}e^{-t^{2}}\ell+O(\log{\ell})
\end{equation}
with some absolute number $c>0$, in contrast to \eqref{eq:nodal length variance Tl}. Moreover, for every $t_{1}, t_{2} \ne 0$, the lengths of the corresponding level curves
asymptically fully correlate ~\cite[Formula (19)]{W Dartmouth}, in the sense that
\begin{equation}
\label{eq:full correlation}
\corr(\len(T_{\ell}^{-1}(t_{1})), \len(T_{\ell}^{-1}(t_{2}))) = 1-o_{\ell\rightarrow\infty}(1),
\end{equation}
where, as usual, for two random, non-degenerate, variables $X,Y$,
$$\corr(X,Y):=\frac{\cov(X,Y)}{\sqrt{\var(X)\cdot\var(Y)}}\in [-1,1].$$

To prove either \eqref{eq:var level set} or \eqref{eq:full correlation} one asymptotically expands
the appropriate the $2$-point correlation function
(resp. cross-correlation function) above Planck scale (see \S\ref{sec:Approx Mixed KR}). No Berry's cancellation in the non-nodal case means that our job here is {\em easier} as compared to \eqref{eq:K2 sphere expand}, as the leading term will give the bulk of contribution
to the variance (resp. covariance), so there is no need for $2$-term expansion.

Finally, asymptotically, as $\ell\rightarrow\infty$, there is no correlation between the nodal length and any level curve at nonzero level $t\ne 0$,
\begin{equation}
\label{eq:nodal vs nonzero}
\corr(\len(T_{\ell}^{-1}(t)), \len(T_{\ell}^{-1}(0))) = o_{\ell\rightarrow\infty}(1).
\end{equation}
While it is easy to check that fact via the cross-correlation route, it is intuitively clear from the corresponding result for
the excursion sets ~\cite{MaWi1}. Because of that it was also later observed ~\cite{MaRoAHP} that the projection of the non-nodal length at level $t\ne 0$ to the $2$nd Wiener chaos, dominating the fluctuations of the level curve length, is, up to
the factor $\sqrt{c}e^{-t^2/2}t^{2}$, equal to the squared norm of $T_{\ell}$. This is in contrast to the nodal case, that is invariant to products by a constant, hence the projection onto the $2$nd Wiener chaos vanishes {\em precisely}, and the fluctuations of the nodal length are dominated by its projection onto the $4$th Wiener chaos. Knowing that, one may also {\em re-prove} ~\cite[Theorem 1.3]{MaRoAHP} both \eqref{eq:full correlation} and \eqref{eq:nodal vs nonzero}.

\subsubsection{Generic nonlinear functionals of $T_{\ell}$ and CLT for nodal length}
\label{sec:gen nonlin func}

Let $H_{q}(\cdot):\R\rightarrow\R$, $q\ge 0$ be the $q$'th degree Hermite polynomials, and consider the centred random variables
\begin{equation*}
h_{q;\ell} := \int\limits_{\Sc^{2}}H_{q}(T_{\ell}(x))dx.
\end{equation*}
For example, $h_{4;\ell}$ is usually referred to as ``sample trispectrum" of $T_{\ell}$ ~\cite{MRW}.
The random variable $h_{q;\ell}$ belongs to the $q$'th Wiener chaos generated by the Gaussian measure induced by $T_{\ell}(\cdot)$.
More generally, let $G:\R\rightarrow\R$ be a function, and, assuming that $$\E[G(T_{\ell}(\cdot))^{2}] <+\infty,$$
define ~\cite{MaWi2} the random variable (``nonlinear functional")
\begin{equation}
\label{eq:Gc gen func}
\Gc(T_{\ell}):= \int\limits_{\Sc^{2}}G(T_{\ell}(x))dx.
\end{equation}

One may Hermite expand
\begin{equation}
\label{eq:G Hermite exp}
G(T_{\ell})\overset{L^{2}(\Sc^{2})}{=\joinrel=}\sum\limits_{q=0}^{\infty}\frac{J_{q}(G)}{q!}H_{q}(T_{\ell}),
\end{equation}
where the coefficients are $$J_{q}(G):=\E[G(Z) \cdot H_{q}(Z)]$$ with $Z$ standard Gaussian, and integrate \eqref{eq:G Hermite exp} on $\Sc^{2}$ to obtain the Wiener chaos expansion
\begin{equation}
\label{eq:Gc gen func sum hq}
\Gc(T_{\ell}) = \sum\limits_{q=0}^{\infty}\frac{J_{q}(G)}{q!}h_{q;\ell}
\end{equation}
for $\Gc(T_{\ell})$.
Hence the $h_{q;\ell}$ are the building blocks of {\em all} square summable functionals of $T_{\ell}$.

If both $q$ and $\ell$ is odd, then, since in this case, $T_{\ell}(\cdot)$ is odd, $h_{q;\ell}\equiv 0$, and, for that reason, in what follows
we assume that $\ell$ is even. It is easy to extend our analysis for either $\ell$ odd, $q$ even, or, otherwise, pass to a subdomain of $\Sc^{2}$ to avoid tautologies
(e.g. restrict $T_{\ell}$ to a hemisphere, as we did above).
It is evident from the properties of Hermite polynomials, that the variance of $h_{q;\ell}$ is, up to an explicit constant, precisely the $q$'th moment of the Legendre
polynomial $P_{\ell}$ on $[0,1]$: $$\var(h_{q;\ell}) =  q!(4\pi)^{2}\int\limits_{0}^{1}P_{\ell}(t)^{q}dt = (4\pi^{2})q!\int\limits_{0}^{\pi/2}P_{\ell}(\cos\theta)^{q}\sin{\theta}d\theta=q!\cdot \int\limits_{\Sc^{2}\times \Sc^{2}}P_{\ell}(\langle x,y\rangle)dxdy,$$
where $\langle \cdot,\cdot\rangle$ is the Euclidean inner product in $\R^{3}$.

\vspace{2mm}

By using Hilb's asymptotics \eqref{eq:Pl Hilb} (and special functions for $q=2,3$), it is possible to evaluate the asymptotics for the moments of the Legendre polynomials, and, therefore, for the variance of $h_{q;\ell}$ as $\ell\rightarrow\infty$ along even integers:
\begin{equation*}
\var(h_{q;\ell}) = \begin{cases} \frac{1}{\ell}(1+o_{\ell\rightarrow\infty}(1)) &q=2 \\
\frac{\log{\ell}}{\ell^{2}}(1+o_{\ell\rightarrow\infty}(1)) &q=4 \\
(4\pi)^{2}q!\frac{c_{q}}{\ell^{2}}(1+o_{\ell\rightarrow\infty}(1)) &q=3,\,q\ge 5
\end{cases}.
\end{equation*}
Here $$c_{q}=\int\limits_{0}^{\infty}\psi J_{0}(\psi)^{q}d\psi\ge 0$$ with the integral conditionally convergent for $q=3$, and absolutely convergent for $q\ge 5$. It is then clear that, for every $q\ge 4$ even,
the a priori nonnegative numbers $c_{q}$ are actually positive, and also ~\cite[Lemma 3.1]{MaWi2} $c_{3} >0$ by evaluating it explicitly via special functions; it is not known to date for which other odd $q$, $c_{q}>0$, though it is believed for all $q\ge 3$.

Further, for every $q\ge 2$ even, for $q=3$, and every $q\ge 5$ so that $c_{q}>0$ a quantitative Central Limit Theorem (q-CLT) was established for $h_{q;\ell}$ as $\ell\rightarrow\infty$ along even $\ell$. Using these results, one may establish the asymptotics of the variance for generic functionals of type \eqref{eq:Gc gen func}, depending on the coefficients of their Hermite expansion \eqref{eq:G Hermite exp}, e.g. if $c_{q}>0$ for at least one $q$
such that $J_{q}(G)\ne 0$. This approach also implies a q-CLT for generic functionals $\Gc(T_{\ell})$ of the said type\footnote{At first, claim these for {\em finite} combinations \eqref{eq:Gc gen func sum hq}.
Then use a standard limit procedure, for sufficiently rapidly coefficients.}.

One in particular important functional of this type is the so-called {\em defect} of $T_{\ell}$, appearing in the physics literature (see e.g. ~\cite{BlGnSm}): It
is the area of the positive excursion set minus the area of the negative excursion set:
\begin{equation*}
\pazocal{D}_{\ell}=\pazocal{D}(T_{\ell}) = \int\limits_{\Sc^{2}}\Hc(T_{\ell}(x))dx,
\end{equation*}
where
\begin{equation*}
\Hc(t) = \begin{cases}
+1 &t>0 \\ -1 &t<0
\end{cases}
\end{equation*}
is the Heaviside function.

The expansion \eqref{eq:Gc gen func sum hq} for the defect only contains odd summands $q\ge 3$, all whose variance is of order of magnitude $\approx \frac{1}{\ell^{2}}$, and there is no single term
dominating the rest of the summands. The asymptotic law of the variance of $ \pazocal{D}_{\ell}$ was previously established to be ~\cite{MaWi3}
\begin{equation*}
\var\left(\pazocal{D}_{\ell}\right) \sim \frac{C}{\ell^{2}},
\end{equation*}
with $C>0$ given by the conditionally converging integral $$C= 32\pi\int\limits_{0}^{\infty}\psi (\arcsin(J_{0}(\psi))-J_{0}(\psi))d\psi  ,$$ incorporating the odd moments $q\ge 3$ of the Legendre polynomials (the Bessel function in the limit). Using the q-CLT for the building blocks\footnote{At least those with $c_{q}>0$; the others have smaller fluctuations, hence, knowing that $c_{3}>0$, these may be neglected.}
$h_{q;\ell}$, a q-CLT for $\pazocal{D}_{\ell}$, via a standard route ~\cite{MaWi2}.

\vspace{2mm}

The functionals \eqref{eq:Gc gen func} do not formally include the nodal length $\len(T_{\ell}^{-1})$ of $T_{\ell}$, since the integral in \eqref{eq:Gc gen func} is only allowed to depend
on the values of $T_{\ell}$, and not its derivatives, cf. \eqref{eq:coarea delta}. Nonetheless, it was observed that $h_{4;\ell}$ asymptotically fully correlates ~\cite[Theorem 1.2]{MRW} with $\len(T_{\ell}^{-1})$ as $\ell\rightarrow\infty$, hence a (quantitative) Central Limit Theorem for $h_{4;\ell}$ will also yield \cite[Corollary 1.3]{MRW} a (quantitative) Central Limit Theorem for $\len(T_{\ell}^{-1})$. This concluded a question opened up a few decades ago.

It is worth mentioning that the projection of $\len(T_{\ell}^{-1}(0))$ to all the odd Wiener chaos spaces vanish, and, in this scenario,
the projection of $\len(T_{\ell}^{-1}(0))$ to the $2$nd chaos vanishes {\em precisely}. The precise vanishing of the projection of the nodal length (and other quantities) is {\em not inherent}
to Berry's cancellation, as in some other situations it will only be {\em majorized} by the projection onto the $4$th Wiener chaos, such as, for example, the same setting of nodal length of $T_{\ell}$, only this time
restricted to some general position subdomain of $\Sc^{2}$ (also see the discussion below on RWM). The precise vanishing of the projection onto the $2$nd chaos is merely a manifestation of the symmetries of $\Sc^{2}$.
This indicates that the asymptotics at infinity of the $2$-point correlation function, as it was originally done in ~\cite{Berry02,W fluctuations}, is a more natural viewpoint on Berry's cancellation phenomenon
than the Wiener chaos expansion language.

\vspace{2mm}

A. Vidotto ~\cite{Vidotto CLT RWM} showed that, as suggested, the above indicated proof yielding the CLT for the nodal length of $T_{\ell}$ also yields the CLT for Berry's RWM restricted to expanding discs $B(R)$ as $R\rightarrow\infty$, proved via a somewhat different argument ~\cite{PucattiEtAl1}, also involving the Wiener chaos expansion.
A. P. Todino ~\cite{Todino local} extended the said results (logarithmic variance, and CLT for nodal length) for nodal length restricted to a spherical length of radius above Planck scale, i.e. for any sequence of radii
$r_{\ell}$, possibly shrinking, so that $r_{\ell}\cdot\ell\rightarrow\infty$. This is a uniform result, optimal w.r.t. the radius, as if $r\ll \frac{1}{\ell}$, then one does not expect the distribution of nodal length to exhibit a limit law of the said type, because of the comparison, after scaling to the RWM, whence the domain is not expanding.

\subsection{Critical points and nodal intersections of independent copies of random spherical harmonics. Higher dimensional spherical harmonics and monochromatic waves.}
\label{sec:crip pnts indep intersect}

Berry proposed\footnote{This model of complex random functions is very different from the models considered by Sodin-Tsirelson (see e.g. ~\cite{Sodin-Tsirelson} in the pure mathematics literature.} ~\cite[Formula (2)]{Berry02} to model the {\em complex} random waves by two independent copies, $u_{1},u_{2}$, of the RWM $u$ in \S\ref{sec:Eucl Riem}, so their a.s. discrete set of zeros is an intersection of the nodal lines of two independent plane waves. As in case of ``real"
plane wave, Berry found that the leading term in the expansion at infinity of the $2$-point correlation of the zeros of the complex plane wave (``nodal number") is purely oscillatory, and does not contribute to the variance. The bulk of the contribution comes from the non-oscillatory part of the $2$nd leading term, and the nodal number variance restricted to $B(R)$ is asymptotic to
$$\var(\#\{x\in \R^{2}: u_{1}(x)=u_{2}(x)=0\}) \sim c_{0}R^{2}\log{R}$$ with some $c_{0}>0$.

This way, we may also consider the number of points in the intersection of the nodal lines of two independent copies of $T_{\ell}$, or, more generally, independently drawn $T_{\ell_{1}}$ and $T_{\ell_{2}}$ with some $\ell_{1}\le \ell_{2}$. Given all the hard work we did to establish \eqref{eq:K2 sphere expand} (and \eqref{eq:nodal length variance Tl}), evaluating the zero density precisely, the $2$-point correlation function at infinity after scaling, and the Wiener chaos expansion for the zeros of $$x\in\Sc^{2}\:\mapsto (T_{\ell_{1}}(x),T_{\ell_{2}}(x)),$$ for the regime $\ell_{1}\rightarrow\infty$, $\ell_{2}=\ell_{2}(\ell_{1})\ge \ell_{1}$, amounts in essence to evaluating the asymptotics of the  determinant of a block diagonal perturbation of the unit matrix, given the determinants of the individual blocks. By this method, it is easy to evaluate the expected number of points in $T_{\ell_{1}}^{-1}(0)\cap T_{\ell_{2}}^{-1}(0) \subseteq\Sc^{2}$ to be proportional to $\ell_{1}\cdot \ell_{2}$ (more precisely, $*\sqrt{\lambda_{\ell_{1}}\cdot \lambda_{\ell_{2}}}$, in line with
\eqref{eq:exp nod length Tl}).

Further, in the said regime, under further assumption that $\ell_{2}$ is not much bigger than $\ell_{1}$, e.g. $\ell_{2}\le C\cdot \ell_{1}$ for some $C>0$ (or beyond it), one may derive the asymptotics
for the variance of the type
\begin{equation*}
\var(\#T_{\ell_{1}}^{-1}(0)\cap T_{\ell_{2}}^{-1}(0)) = c_{1}(\ell_{1}^{2}+\ell_{2}^{2})\log\ell_{1} + c_{2}\ell_{1}\ell_{2}\log\left(\frac{\ell_{1}}{\ell_{2}+1-\ell_{1}}\right),
\end{equation*}
for some explicitly computable constants $c_{1},c_{2}>0$; in particular,
\begin{equation*}
\var(\#T_{1;\ell}^{-1}(0)\cap T_{2;\ell}^{-1}(0)) \sim c_{0}\ell^{2}\log{\ell};
\end{equation*}
the purely oscillatory nature of the leading term for $K_{2}$ indicates that Berry's cancellation occurs in this case again.
If, on the other hand, $\ell_{2}$ is bigger by order of magnitude than $\ell_{1}$, then, at some point, a {\em phase transition}
occurs, which is easy to explicate by analyzing the resulting Kac-Rice integral.
To derive an asymptotic law for $\#T_{\ell_{1}}^{-1}(0)\cap T_{\ell_{2}}^{-1}(0)$ one uses a Wiener chaos expansion, that, in light of the above, is, in the said regime, most certainly dominated by the $4$th chaos projection, analysis of which will yield the Central Limit Theorem for the said quantity.

What is more technically demanding than the analysis of intersection of the nodal lines of random independent copies of $T_{\ell}$ is the analysis of the number of critical points and critical values (i.e. the values of $T_{\ell}(x_{k})$ where $\{x_{k}\}\subseteq\Sc^{2}$ are the critical points of $T_{\ell}$) of $T_{\ell}$, since this amounts to the intersection of nodal lines of two {\em dependent} random Gaussian functions, namely, the derivatives of $T_{\ell}$, so that the associate covariance matrix is no longer block diagonal. These dependencies cannot be neglected, and do impact the asymptotic distribution of the number of critical points (and their values). It is easy to evaluate the expected total number of critical points to be proportional to $\lambda_{\ell}\sim \ell^{2}$, and, by expanding the $2$-point correlation function at infinity, their variance is proportional ~\cite{CW} to $\ell^{2}\log{\ell}$, with Berry's cancellation. For critical values on $\R$, a non-universal limit distribution has been derived ~\cite{CMW} as a {\em point-process}, and a variance for the
number of critical values belonging to a generic interval $I=[a,b]\subseteq\R$ has been established ~\cite{CMW} to be proportional to $\ell^{3}$; for non-generic $I$ the variance is of lower order, most likely (this special class of $I$ is reflecting their symmetry, for example, as it was mentioned, $I=\R$ has a lower order variance of $\ell^{2}\log\ell$, as it is expected from other intervals in this class).
Attraction and repulsion of critical points (also finer, maxima, minima and saddles) were studied for Berry's RWM ~\cite{BCW1}, and generic Gaussian isotropic processes ~\cite{BCW2,AzDel}.

\subsection{Generalizations}

All the above results generalize for $n$-dimensional sphere, $n\ge 2$ and any number $k\le n$ of independent spherical harmonics of either equal or non-equal degrees,
or monochromatic waves on $\R^{n}$, see e.g. ~\cite{PucattiEtAl1,EstradeEtAl}, whence one is interested in the corresponding $(n-k)$-volume (number of discrete zeros for $k=n$);
every configuration $(n,k)$ is susceptible to Berry's cancellation. Alternatively, one may consider
intersections of level sets: The analogous results to \S\ref{sec:nodal volume} (especially towards the end) will hold in this more general setting, with the only outstanding question, being somewhat vague, is that of the number of points
$\# T_{1;\ell}^{-1}(0)\cap T_{2;\ell}^{-1}(t)$ in the intersection
of the nodal line of $T_{\ell}$ and level $t\ne 0$ curve of its independent copy. Without performing the necessary (straightforward) computation, it seems that this case will exhibit some features of nonzero level curve length of a single spherical harmonic and than that of the nodal length (in terms of Berry's cancellation and the dominating chaos projection).

\subsection{Nodal intersections against smooth curves}

\label{sec:nod intersect RWM}

Let $\gamma:[0,L]\rightarrow\R^{2}$ be a (fixed) smooth curve\footnote{It seems that $C^{2}$-smoothness should be sufficient},
arc length parameterized. In this context, one is interested in the number of the intersections of Berry's RWM $u(\cdot)$ with
$\gamma$, or, rather, of the wavenumber-$k$ plane waves $u(k\cdot)$, $k\rightarrow\infty$. Equivalently, this is the number of the a.s. discrete zeros of the function $$f_{k}(t):=u(k\gamma(t)):[0,L]\rightarrow\R.$$ Note that the random Gaussian process is not stationary. However, it is still quite easy to evaluate the zero density $K_{1}(t)=K_{1;f_{k}}(t)$ in this case, and exploit the univariate nature of $f_{k}$ to find that
it is universally given by
\begin{equation*}
K_{1}(t) \equiv \frac{k}{\sqrt{2}\pi},
\end{equation*}
independent of the geometry of $\gamma$, so that the expected number of nodal intersections in this case is given precisely by
\begin{equation*}
\E[\#f_{k}^{-1}(0)] = \frac{k}{\sqrt{2}\pi}\cdot L.
\end{equation*}

Evaluating the variance of the nodal intersections number is much more technically demanding and subtle. Nevertheless it was successfully performed to show that, as $k\rightarrow\infty$,
\begin{equation*}
\var(\#f_{k}^{-1}(0)) \sim \frac{L}{2\pi^{3}}k\log{k}.
\end{equation*}
A few observations are due. Somewhat surprisingly, the leading asymptotics of the variance is universally independent of the geometry of $\gamma$. This is especially unexpected, or even shocking, due to the existence of the so-called ``static" curves in case of Arithmetic Random Waves (see \S\ref{sec:ARW ni} below) that exhibit lower variance than the generic curves, and exhibit non-central limit law. Despite the existence of logarithmic term in
the variance, there is no Berry's cancellation here, and the bulk of the contribution to the variance arises from the leading term of the expansion of $K_{2}(t_{1},t_{2})$ away from the diagonal, i.e. $|t_{2}-t_{1}|\gg 1/k$. This indicates that, most certainly, the projection to the $2$nd Wiener chaos of $\#f_{k}^{-1}(0)$ will dominate its fluctuations, and will yield the quantitative CLT via the standard route. It is highly likely that same techniques yield the {\em small scale} analogues of these all the way to the Planck scale, cf. ~\cite[Theorem 1.2]{GWAJM}. Finally, same results are, with no doubt,
also applicable for the nodal intersections of the random spherical harmonics $T_{\ell}$ against a fixed smooth spherical curves. Higher dimensional analogues should also be treatable via similar methods.

\subsection{Berry's cancellation - generic Gaussian random fields on Euclidean space}

Let $F$ be a stationary random Gaussian random field on $\R^{2}$, and let $\rho$ be its spectral measure, i.e. $\rho$ is the Fourier transform of the covariance function of $F$.
Assuming that $F$ is {\em monochromatic}, i.e. $\rho$ is supported on $\mathcal{S}^{1}\subseteq\R^{2}$, the covariance function of $F$ is given by the inverse Fourier transform
\begin{equation}
\label{eq:rF Fourier transform circle}
r_{F}(x) = \int\limits_{\Sc^{1}}e^{i\langle x,\vartheta\rangle}d\rho(\vartheta),
\end{equation}
where $x=(x_{1},x_{2})$, $\vartheta=(\cos\theta,\sin\theta)\in\Sc^{1}$, $\theta\in [0,2\pi]$, $\langle x,\vartheta\rangle = x_{1}\cos\theta+x_{2}\sin\theta$, and we may further assume w.l.o.g.
that $\rho$ is invariant w.r.t. rotation by $\pi$. It is then straightforward to obtain a precise expression for the expected nodal length of $F$ restricted to $B(R)$ to be $c(\rho)\cdot \pi R^{2}$,
where $c(\rho)$ is explicitly computable in terms of a couple of power moments of $\rho$. The question put forward is, again, that of the asymptotic variance (and the asymptotic law),
as $R\rightarrow\infty$, the intrigue being whether Berry's cancellation occurs, and in what generality if it does.

One may obtain an asymptotic expansion of $r_{F}(x)$ and its derivatives at infinity, i.e. when $|x|\rightarrow\infty$, in analogy to \eqref{eq:J0 2 term asymp}
for the isotropic case \eqref{eq:RWM covar} of Berry's RWM, and use these to obtain an asymptotic expression at infinity for $K_{F;2}(\cdot)$, the $2$-point correlation function \eqref{eq:K2 def Gauss int}.
First, we assume that $\rho$ has a density on $\Sc^{1}$: $d\rho(\theta)=g(\theta)d\theta$ for some nonnegative ``nice" function $g:[0,2\pi)\rightarrow\R_{\ge 0}$ ($C^{2}$ should be sufficient).
In this case, it is possible to use the {\em stationary phase} method on \eqref{eq:rF Fourier transform circle} to find that, after cancelling out the squared expectation, the leading term in the expansion of $K_{F;2}(x)$ at infinity is of purely oscillatory nature, and, using the polar representation $x=|x|e^{i\phi}$, the next term is proportional to $$g(\phi)^{4}\cdot \frac{1}{\|x\|^{2}}.$$ Hence this model is susceptible to Berry's cancellation. By integrating \eqref{eq:KC 2nd mom} the latter asymptotics, and subtracting the squared expectation, we may then deduce that, as $R\rightarrow\infty$,
\begin{equation*}
\var(F^{-1}(0)\cap B(R)) \sim c(g)\cdot R^{2}\log{R},
\end{equation*}
with
\begin{equation}
\label{eq:c(g) const 4th mom}
c(g)=c_{0}\cdot \int\limits_{\Sc^{1}}g(\theta)^{4}d\theta,
\end{equation}
and $c_{0}>0$ an absolute constant. The leading constant \eqref{eq:c(g) const 4th mom} is reminiscent of the previously discovered one for the Arithmetic Random Waves (see \eqref{eq:var ARW KKW} in \S\ref{sec:survey ARW} below), although it bears the $4$th Fourier component of the spectral measure, as opposed to the $4$th moment \eqref{eq:c(g) const 4th mom} of its density. Note that, interestingly, the expression
$g(\theta)^{4}$ is minimized by the {\em uniform} probability measure, due to Jensen's inequality; somewhat surprisingly, the leading constant \eqref{eq:c(g) const 4th mom} does not bear the derivatives of $g$,
despite the contrary for the intermediate expressions (various cancellations yielding a very neat and compact expression for the leading non-oscillatory asymptotics of the $2$-point correlation, and the nodal length variance).

There is no doubt that, using the Wiener chaos expansion for $\var(F^{-1}(0)\cap B(R))$ reveals that, in this case, the projection onto the $4$th Wiener chaos is dominating the projection onto the $2$nd one (which asymptotically vanishes, due to the oscillatory nature of the $2$ leading term of $K_{F;2}$), and also the higher Wiener chaos spaces. A Central Limit Theorem for $\var(F^{-1}(0)\cap B(R))$ also follows from the Central Limit Theorem for the projection onto the $4$th Wiener chaos, via the standard route indicated above. If, on the other hand, $\rho$ has jump discontinuities, or $\rho$ is either purely atomic, or $\rho$ is a superposition of atoms and an absolutely continuous part, then the atoms would dominate the asymptotics at infinity of the covariance function \eqref{eq:rF Fourier transform circle}, and the nodal length variance, though some notion of Berry's cancellation will be holding in some generality (for example, in case of jump discontinuities only, the variance will be logarithmic, with leading constant depending on the jumps). Since if the support of $\rho$ has a point
of interior, $r_{F}(\cdot)$ and its derivatives have rapid decay in infinity, the nodal length variance will be of smaller order of magnitude $R^{2}$, with no Berry's cancellation.

In this light we may formulate the following {\em principle} or {\em meta-theorem}: {\bf Under appropriate assumptions on $F$ (or $\rho$), it is susceptible to Berry's cancellation, if and only if $F$ is monochromatic.}
All of the above (including the results stated in previous sections) could be generalized to any dimensions $n\ge 2$, and intersections of $k\le n$ nodal sets of independent copies of $F$ (or different functions), or other related quantities associated with $F$, such as, for example the number of critical points.

\vspace{2mm}

Finally, recall that the ensemble $T_{\ell}(\cdot)$ scales like Berry's RWM around {\em every} reference point $x\in\Sc^{2}$.
It would be interesting, given a monochromatic random field $F:\R^{2}\rightarrow\R$ of the said type, to construct an ensemble of Gaussian functions $T_{F;\ell}(\cdot)$, supported on
the space of degree $\ell$ spherical harmonics, that scales like $F$ around every (or, at least, almost every) reference point in $\Sc^{2}$.

\subsection{Band-limited functions - general case}

\label{sec:band-lim survey}

Let $\M$ be a smooth closed $n$-manifold with $n\ge 2$, i.e. $\M$ is a smooth compact Riemannian manifold with empty boundary. Under further assumption\footnote{This condition is satisfied for ``most" of the manifolds in
many different senses.}
that the measure of the geodesic loops through a point $x\in\M$ is $0$ for a.a. $w\in\M$, let
\begin{equation*}
f_{\lambda}(x) = \sum\limits_{\lambda_{i}\in [\lambda,\lambda+1]}c_{i}\varphi_{i}(x),
\end{equation*}
with the $c_{i}$ i.i.d. standard Gaussian r.v., i.e. $f_{\lambda}(\cdot)$ are the random Gaussian monochromatic waves
on $\M$, corresponding to the shortest possible energy window. Alternatively, without the geometric restriction, take
\begin{equation*}
f_{\lambda}(x) = \sum\limits_{\lambda_{i}\in [\lambda-\eta(\lambda),\lambda]}c_{i}\varphi_{i}(x),
\end{equation*}
with some $\eta\rightarrow\infty$ (but $\eta=o(\lambda)$).

Zelditch ~\cite{Zelditch Sunada} showed that the total expected nodal $(n-1)$-volume of $f_{\lambda}$ is asymptotic to
\begin{equation*}
\E[\vol_{n-1}(f_{\lambda}^{-1}(0))] = c_{\M} \cdot \lambda,
\end{equation*}
with $c_{\M}>0$ explicit, proportional to the Riemannian volume of $\M$.
Canzani-Hanin ~\cite{CH} employed the $2$-point correlation function ~\eqref{eq:K2 def Gauss int} to give an upper bound for the variance
$\var(\vol_{n-1}(f_{\lambda}^{-1}(0)))$. They developed a sophisticated partition argument, carefully bounding the contribution of the $K_{2}(\cdot,\cdot)$ in several ranges,
to yield the handy bound ~\cite[Theorem 1]{CH}
\begin{equation}
\label{eq:CH upper bnd}
\var(\vol_{n-1}(f_{\lambda}^{-1}(0))) = O(\lambda^{2-(n-1)/2}),
\end{equation}
or, equivalently, to show that the {\em normalized} variance
\begin{equation*}
\var\left(\frac{\vol_{n-1}(f_{\lambda}^{-1}(0))}{\lambda}\right)= O\left( \lambda^{-(n-1)/2} \right)
\end{equation*}
vanishes.

It seems that, by employing the {\em Approximate} Kac-Rice technique as in \S\ref{sec:Approx Mixed KR} in place of the (precise) Kac-Rice, one can simplify or shorten Canzani-Hanin's argument validating
the sufficient Gaussian non-degeneracy conditions for the Kac-Rice to hold. Instead, it is sufficient to prove the same for {\em most} $(x,y)\in \M\times \M$. Further, by writing
the Kac-Rice formula and using the upper bound in terms of the $L^{2}$-energy of the covariance function and its derivatives of the type \eqref{eq:var bnd L2 enrgy covar}, it is expected to yield a sharper upper bound
than \eqref{eq:CH upper bnd} in terms of the power of $\lambda$, not unlikely
\begin{equation*}
\var(\vol_{n-1}(f_{\lambda}^{-1}(0))) = O(\lambda).
\end{equation*}
Unfortunately, all of these shed no light on the true behaviour of the nodal volume variance, nor give a nontrivial {\em lower} bound for the said variance, ditto of the same order of magnitude as the upper bound, cf. \eqref{eq:nodal length variance Tl}. That the {\em local} convergence of the covariance (even on a logarithmic scale above the Planck scale $r\approx 1/\lambda$) and its derivatives to those of the RWM does not indicate the same for the total nodal length follows from the toral application, see \S\ref{sec:nod length}.

\vspace{2mm}

Let $x_{0}\in \M$ be a reference point, and consider the function $f_{\lambda}$ restricted to the geodesic ball $B_{x_{0}}(1/\lambda)$. It follows that the scaled random field \eqref{eq:band-lim scale g},
defined on $B(R)$, converges to the ($n$-dimensional) RWM $g_{\infty}$. One may then use the Continuous Mapping Theorem ~\cite{CH} (see also the proof in ~\cite[p. 70]{Billingsley}),
to infer that, up to the scaling, the distribution of the nodal volume of $f_{\lambda}$ restricted to $B(R/\lambda)$ converges in law to the distribution of the nodal volume of $g_{\infty}$ on $B(R)$ (with no information what this distribution is). That the sufficient conditions
for the application of the Continuous Mapping Theorem are satisfied for the nodal volume of $g_{\lambda}$ (i.e. the nodal volume is a.s. continuous w.r.t. $g_{\lambda}$) follows from an application of the well-known Bulinskaya's Lemma. Since the said result holds for every $R>0$, one may also {\em tour de force} the Continuous Mapping Theorem while using Berry's asymptotics \eqref{eq:Berry log} and CLT for the nodal length variance of the RWM in $2$d, to infer the analogous result for the nodal length of $f_{\lambda}$ restricted to geodesic balls in $\M$ of radius $R/\lambda$, with $R$ slowly growing (with no control how slow the growth of $R$ should be).

Finally, by using an extension due to Keeler ~\cite{Keeler} of the scaling property \eqref{eq:covar band-lim->covar infty} for $\alpha=1$,
holding for a further restricted (but still generic) class of manifold, a logarithm power above the Planck scale, one may pushforward
~\cite{PucattiEtAl2} the asymptotic variance and the CLT logarithmically above the Planck scale, for the said class of manifolds. As the Arithmetic Random Waves of \S\ref{sec:survey ARW} is an ensemble of random Gaussian band-limited functions that locally scales like the RWM for a generic sequence of toral energies, though the asymptotic distribution of the total nodal length is very different, we conclude that the above is a local phenomenon, generally not indicative of the global nature of the nodal line.

\section{Survey of results: Arithmetic random waves}
\label{sec:survey ARW}

\subsection{Lattice points: angular distribution, spectral correlations, quasi-correlations, and semi-correlations}

Recall that $f_{n}$ are the Arithmetic Random Waves (ARW) given by \eqref{eq:fn ARW def}, with covariance function given by \eqref{eq:rn ARW covar}. Unlike the spherical harmonics, unfortunately, $r_{n}(\cdot)$ (resp. the mixed derivatives of $r_{n}$) does
not admit an asymptotic law (after scaling), valid on the whole of $\Tb^{2}$, cf. \eqref{eq:covar Legendre} and \eqref{eq:Pl Hilb}. As a substitute, one can work ~\cite{KKW} with the {\em moments} of $r_{n}(\cdot)$ (resp. its derivatives) on $\Tb^{2}$. For an integer number $k\ge 2$ one has
\begin{equation}
\label{eq:moments rn correlations}
\int\limits_{\Tb^{2}}r_{n}(x)^{k}dx =  \frac{1}{N_{n}^{k}}\#\Rc_{k}(n),
\end{equation}
where
\begin{equation}
\label{eq:Rc correlations}
\#\Rc_{k}(n) := \{(\lambda_{1},\ldots,\lambda_{k})\in\Lambda_{n}^{k}:\lambda_{1}+\ldots+\lambda_{k}=0 \}
\end{equation}
is the length-$k$ {\em spectral correlation set}.

The sets $\Rc_{k}(n)$ (equivalently, the moments of $r_{n}(\cdot)$) play a crucial role in the analysis of the nodal length of $f_{n}$.
For $k=2$, trivially $$\Rc_{2}(n)=\{(\lambda,-\lambda):\: \lambda\in\Lambda_{n}\}.$$ The analogue of this holds as well for $k=4$, thanks to Zygmund's trick ~\cite{Zygmund} observing that given $\lambda_{1},\lambda_{2}\in\Lambda_{n}$,
the $\lambda_{3},\lambda_{4}\in\Lambda_{n}$ are determined, up to a sign and a permutation. That is,
the length-$4$ spectral correlation set equals to the diagonal set, i.e. it is the set $$\Rc_{4}(n)=\{\pi(\lambda,-\lambda,\lambda',-\lambda'):\: \lambda,\lambda'\in\Lambda_{n},\pi\in S_{4}\}$$ of permutations of those $4$-tuples cancelling out in pairs, hence its cardinality is
\begin{equation*}
\label{eq:Zygmund}
\#\Rc_{4}(n) = 3N_{n}(N_{n}-1).
\end{equation*}

For our purposes we need a nontrivial bound for $\#\Rc_{6}(n)$. As above, given $\lambda_{1},\ldots,\lambda_{4}\in\Lambda_{n}$, it determines, up to a sign and a permutation, $\lambda_{5},\lambda_{6}\in\Lambda_{n}$. Therefore, it follows that $\#\Rc_{4}(n)=O(N_{n}^{4})$. For our needs, any improvement of the type
\begin{equation}
\label{eq:Rc6=o(N^4)}
\#\Rc_{6}(n) = o_{N_{n\rightarrow\infty}}(N_{n}^{4})
\end{equation}
is required. The result \eqref{eq:Rc6=o(N^4)} was established by Bourgain, published as part of ~\cite{KKW}. Since then, a better bound of
\begin{equation*}
\#\Rc_{6}(n) = O(N_{n}^{7/2})
\end{equation*}
was proven ~\cite{BB}, using the power of the ring of Gaussian integers $\Z[i]\cong \Z^{2}$. Some better bounds of
\begin{equation*}
\#\Rc_{6}(n) = O(N_{n}^{3})
\end{equation*}
were also asserted ~\cite{BB}, albeit either conditional, or valid unconditionally for a generic subsequence of $S$.

\vspace{2mm}

Let us say a few words about the dimension $N_{n}=r_{2}(n)$ of the eigenspaces of toral eigenfunctions. It fluctuates in an erratic way. First, it is easy to show that $N_{n}=O(n^{o(1)})$, i.e., that for every $\epsilon>0$, $N_{n}=O(n^{\epsilon})$. On the one hand, $N_{n}$ could grow as a power of $\log{n}$, and its normal order is $N_{n}\approx \log{n}^{\log{2}/2}$, that is, for every $\epsilon>0$,
$$ \log{n}^{\log{2}/2-\epsilon} \le N_{n} \le \log{n}^{\log{2}/2+\epsilon}$$ for a density-$1$ sequence $\{n\}\subseteq S$. On the other hand, $N_{n}$ is as small as $N_{n}=8$ for an infinite sequence of primes
$p\equiv 1 \mod 4$. We do assume throughout that $N_{n}\rightarrow\infty$.

\vspace{2mm}

To measure the angular distribution $$\{\lambda/\sqrt{n}:\: \lambda\in\Lambda_{n}\}\subseteq \Sc^{1}$$ of the lattice points $\Lambda_{n}$ lying on the circle $\sqrt{n}\Sc^{1}$, one introduces ~\cite{KKW} the atomic probability measure on $\Sc^{1}$
\begin{equation}
\label{eq:mun atom measure}
\mu_{n}:=\frac{1}{N_{n}}\sum\limits_{\lambda\in\Lambda_{n}}\delta_{\lambda/\sqrt{n}}
\end{equation}
supported on the angles corresponding to these lattice points. It is well-known ~\cite{KK,EH} that, for a density-$1$ sequence, the angles of $\Lambda_{n}$ become {\em equidistributed} in a quantitative sense, so that, in particular,
\begin{equation}
\label{eq:equidistribution KK,EH}
\mu_{n}\Rightarrow \frac{d\vartheta}{2\pi}.
\end{equation}
This also means that the covariance function, $r_{n}$ as in \eqref{eq:rn ARW covar}, of the ARW, tends, {\em locally at Planck scale}, to the covariance function \eqref{eq:RWM covar} of Berry's RWM, together with their respective derivatives. Finer than that, the {\em quantitative} equidistribution asserts that the same is true at radii a power of logarithm {\em above} Planck scale.

On the other hand, other limit measures in \eqref{eq:equidistribution KK,EH} are possible: Cilleruelo ~\cite{Cilleruelo} constructed a ``Cilleruelo" sequence $\{n\}\subseteq S$ so that
\begin{equation*}
\mu_{n}\Rightarrow \frac{1}{4}\left(\delta_{\pm 1}+\delta_{\pm i}\right),
\end{equation*}
where we think of $\Sc^{1}\subseteq\C$. The rich structure of all possible measure ``attainable" as a limit of a subsequence of $\mu_{n}$ was studied ~\cite{KuWi,SartoriAttainable}.

\vspace{2mm}

Another subject, naturally occurring when studying the eigenfunctions on a small disc of radius above Planck scale is that of {\em quasi-correlations}. For $k\ge 2$, and $\epsilon>0$,
define the quasi-correlation set as
\begin{equation*}
\Qc_{k}(n,n^{1/2-\epsilon}) = \left\{(\lambda_{1},\ldots,\lambda_{k})\in\Lambda_{n}^{k}:\: 0<\left\|\sum\limits_{j=1}^{k}\lambda_{j}\right\| < n^{1/2-\epsilon}\right\},
\end{equation*}
i.e. those $k$-tuples whose sum has cancellations of order of magnitude $n^{\epsilon}$. For example, for $k=2$, the quasi-correlations are identified as the close-by pairs, and it was asserted ~\cite{GWCMP} that for a generic $n\in S$, $$\Qc_{2}(n,n^{1/2-\epsilon})=\varnothing,$$ and the same was inferred ~\cite[Lemma 5.2]{RoWi} for $k=4$. Additionally, the notion of {\em semi-correlations} was introduced
~\cite{CamKlurWig}\footnote{This is an opportunity to correct an informal statement made in that manuscript. There are two negative biases of different nature. The perpendicular intersection of the nodal line with the boundary has only a local effect on the density in its neighbourhood, and is unrelated to the negative excess term infinitely many wavelengths away from the boundary.}, i.e. those $k$-tuples $(\lambda_{1},\ldots,\lambda_{k})\in\Lambda_{n}^{k}$, with the $1$st {\em coordinate} of $\lambda_{1}+\ldots+\lambda_{k}$ vanishing, for the study of the {\em boundary-adapted} ARW, i.e. random Laplace eigenfunctions on the unit square.

\subsection{Nodal length of Arithmetic Random Waves}
\label{sec:nod length}

It is easy to compute the expected nodal length of $f_{n}$ by invoking the stationarity of $f_{n}$, so that the zero density $K_{1}(\cdot)$ depends only on $n$. Hence,
an explicit computation of the involved Gaussian integral yields a precise expression for the expected nodal length for every $n\in S$ to be ~\cite{RudWig1}
\begin{equation*}
\E[\len(f_{n}^{-1}(0))] = \frac{\pi}{\sqrt{2}}\sqrt{n}.
\end{equation*}
A handy upper bound was also given for the variance
\begin{equation*}
\var\left(\len(f_{n}^{-1}(0))\right) = O\left(\frac{n}{\sqrt{N_{n}}}   \right),
\end{equation*}
as $N_{n}\rightarrow\infty$. A precise asymptotic law ~\cite{KKW} of the form
\begin{equation}
\label{eq:var ARW KKW}
\var\left(\len(f_{n}^{-1}(0))\right) = (1+\widehat{\mu_{n}}(4)^{2})\frac{\pi^{2}}{128}\frac{n}{N_{n}^{2}}(1+o_{N_{n}\rightarrow\infty}(1))
\end{equation}
was finally derived for the variance. To control the error term in the asymptotics \eqref{eq:var ARW KKW}, a crucial role has been played by a nontrivial upper bound
\begin{equation*}
\#\Rc_{6}(n) = o_{N_{n}\rightarrow\infty}(N_{n}^{4})
\end{equation*}
due to J. Bourgain.

A few observations on \eqref{eq:var ARW KKW} are due. Remarkably, \eqref{eq:var ARW KKW} shows that, assuming
$N_{n}\rightarrow\infty$, the corresponding sequence of variances $\var\left(\len(f_{n}^{-1}(0))\right)$ {\em fluctuates} depending on the angular distribution of the lattice points $\Lambda_{n}$. Thus,
to exhibit an asymptotic law, it is essential to pass to a subsequence of $S$ so that the corresponding $\Lambda_{n}$ would observe a limit angular distribution, i.e. $\mu_{n}\Rightarrow\mu$ for some
probability measure $\mu$ on $\Sc^{1}$. The variance \eqref{eq:var ARW KKW} is of a smaller order of magnitude as compared to what was the expected order of magnitude ~\cite{RudWig1}
of $\frac{n}{N_{n}}$. This is due to another unexpected cancellation (``arithmetic Berry's cancellation"), of a similar nature, but different appearance, to Berry's cancellation \eqref{eq:Berry log}
in the context of RWM, or random spherical harmonics. Similar to random spherical harmonics,
Berry's cancellation does not occur for level curves at nonzero level: For every $t\ne 0$,
\begin{equation*}
\var\left(\len(f_{n}^{-1}(t))\right) \sim c_{0}t^{4}e^{-t}\frac{n}{N_{n}}
\end{equation*}
as $N_{n}\rightarrow\infty$, with some absolute constant $c_{0}>0$ (cf. \eqref{eq:var level set}). Further, the full correlation phenomenon \eqref{eq:full correlation} is observed
in this context too:
for every $t_{1},t_{2} \ne 0$:
\begin{equation*}
\corr(\len(f_{n}^{-1}(t_{1})), \len(f_{n}^{-1}(t_{2}))) = 1-o_{N_{n}\rightarrow\infty}(1).
\end{equation*}

\vspace{2mm}

To study the fluctuations of $\len(f_{n}^{-1}(0))$ in more details, the Wiener chaos expansion technique was applied. In accordance to the analogous situation for the spherical harmonics
(see \S\ref{sec:gen nonlin func}), the $2$nd Wiener chaos component of $\len(f_{n}^{-1}(0))$ vanishes precisely, whence the bulk of the contribution to the fluctuations of $\len(f_{n}^{-1}(0))$ comes from its projection to the $4$th Wiener chaos component; bounding the contribution of the higher Wiener chaos spaces amounts to bounding the size of the spectral correlations $\Rc_{k}(n)$ in \eqref{eq:Rc correlations}, $k\ge 6$, $n\rightarrow\infty$, an approach that could also be used re-prove the variance estimate \eqref{eq:var ARW KKW}, implemented within \cite{PRAbel}.
A non-universal non-Central Limit theorem was established ~\cite{MRWGafa} for
$$ \widetilde{\Lc_{n}}:=\frac{\len(f_{n}^{-1}(0)) - \E[\len(f_{n}^{-1}(0))]}{\sqrt{\len(f_{n}^{-1}(0))}}$$ for the regime $N_{n}\rightarrow\infty$, where to exhibit a limit law, it is essential to separate the sequence
$S$ according to the angular distribution of the lattice points $\Lambda_{n}$, in the following, explicit, manner. For $\eta\in [-1,1]$ denote the random variable
$$M_{\eta}:= \frac{1}{2\sqrt{1+\eta^{2}}}\left( 2-(1+\eta)X_{1}^{2}-(1-\eta)X_{2}^{2} \right),$$ where $X_{1},X_{2}$ are i.i.d. standard Gaussian r.v. Then, for $\eta\in [-1,1]$, if $\{n\}\subseteq S$ is a subsequence of energies so that $N_{n}\rightarrow\infty$ and
$$\widehat{\mu_{n}}(4)\rightarrow \eta,$$ one has
\begin{equation}
\label{eq:Lc(n)->Meta}
\widetilde{\Lc_{n}} \rightarrow M_{\eta},
\end{equation}
convergence in distribution (and also in some other senses).

\vspace{2mm}

As the variance and the limit law of the {\em total} nodal length of $f_{n}$ are completely understood, the natural question is that of the {\em local} analogues of these, above Plack scale, i.e. the nodal length of $f_{n}$ restricted to a disc of radius $r=r(n)$ so that $r_{n}\cdot \sqrt{n}\rightarrow\infty$; by the stationarity of $f_{n}(\cdot)$, one may only consider the {\em centred} discs $B(r_{n})=B_{0}(r_{n})$. As in the total length case, it is also naturally invoking the moments of the covariance $r_{n}$ in \eqref{eq:rn ARW covar} and its derivatives, this time {\em restricted} to $B(r_{n})$ (or, rather, to $B(r_{n})\times B(r_{n})$). Fix a small number $\epsilon>0$, and consider $r_{n}:=n^{-1/2+\epsilon}$. The analogue equality to
\eqref{eq:moments rn correlations} holds, this time, up to an error term that is expressed in terms of the size of the {\em quasi-correlations} set, with $\delta=n^{1/2-\epsilon}$. It was shown ~\cite[Theorem 1.4]{BMW} that, for a ``generic" ({\em density}-$1$) sequence, for given a number $k\ge 2$, the spectral quasi-correlation set of length $m=2,\,\ldots,\, k$ is {\em empty}, and, in addition, such a generic sequence $\{n\}$ could be selected sufficiently rich for the accumulation set of $\{\widehat{\mu_{n}}(4)\}$ to contain the interval $[0,1]$.

Let $\Lc_{n}$ be the total nodal length of $f_{n}$, $\epsilon>0$, $r_{n}=n^{-1/2+\epsilon}$, and $\Lc_{n,r_{n}}$ be the nodal length of $f_{n}$ restricted to $B(r_{n})$. It was shown that if $\{n\}\subseteq S$ is a sequence of energies, with the quasi-correlation sets of both lengths $2$ and $6$ sufficiently small\footnote{1. This, in particular, is satisfied when the quasi-correlations sets are empty. 2. In particular, this implies $N_{n}\rightarrow\infty$.}, then
\begin{equation}
\label{eq:corr loc tot nodal len ARW}
\corr(\Lc_{n},\Lc_{n,r_{n}}) = 1-o_{n\rightarrow\infty}(1).
\end{equation}
In particular, the local analogues of \eqref{eq:var ARW KKW} and \eqref{eq:Lc(n)->Meta} are satisfied just above the Planck scale, for a generic sequence of energy levels, with $\eta$ (rather, $\eta^{2}$) unrestricted.
The full correlation \eqref{eq:corr loc tot nodal len ARW} was proved by using an appropriate {\em Mixed} Kac-Rice formula, to evaluate the relevant covariance, and upon using \eqref{eq:var ARW KKW}, and its local analogue that is the bulk of the proof. One reuses the result on the $2$-point correlation function that was instrumental for establishing \eqref{eq:var ARW KKW}, to make a $2$-line proof of \eqref{eq:corr loc tot nodal len ARW}
given all the developed technology. The alternative for demonstrating this result is to use the Wiener chaos expansion approach, that, up to an admissible error and explicit factors, identifies the Wiener chaos expansion of
$\Lc_{n,r_{n}}$ with that of $\Lc_{n}$. This also shows that the bulk of the fluctuations of $\Lc_{n,r_{n}}$ is its projection onto the $4$th Wiener chaos space, though, unlike $\Lc_{n}$, its projection onto the $2$nd Wiener
chaos does {\em not} vanish precisely (consistent to the sphere).

Sartori ~\cite[Theorem 1.3]{Sartori} refined the latter result for \eqref{eq:corr loc tot nodal len ARW} (and thus the local variance and limit law) to hold a logarithm power above the Planck scale, i.e. for $r_{n}=\frac{\log^{A}}{n^{1/2}}$
with $A>0$ sufficiently large, instead of a power of $n$ above the Planck scale. On the other hand, the angular equidisribution result \eqref{eq:equidistribution KK,EH} holding for generic $n$ means that,
for a density-$1$ sequence $\{n\}\subseteq S$, the restriction $f_{n}|_{B(r_{n})}$ to $B(r_{n})$ the random field $f_{n}|_{B(r_{n})}$ with $r_{n}=\frac{\log^{\delta}}{n^{1/2}}$, $\delta>0$ sufficiently small,
converge to Berry's RWM.
By constructing an appropriate coupling of $f_{n}$ with Berry's RWM, it is therefore possible ~\cite[Theorem 1.4]{Sartori} to infer that, for this choice of $r_{n}$, $\Lc_{n,r_{n}}$ has Berry's logarithmic variance and converges to the Gaussian distribution (see also ~\cite{PucattiEtAl2}). Therefore,
the Arithmetic Random Waves is one example, when the local distribution of the nodal length is not indicative of the total nodal length (see the discussion at the end of \S\ref{sec:band-lim survey}).

\vspace{2mm}

Benatar and Maffucci ~\cite{BeMa} considered the $3$-dimensional Arithmetic Random Waves, i.e. the $3$-dimensional random toral eigenfunctions.
On one hand, for the $3$d torus, the angles of the lattice points lying on spheres of growing radii equidistribute ~\cite{Golubeva,Duke} for the {\em full} sequence of energy levels (i.e. the numbers representable as sum of $3$ squares), e.g. in the sense of convergence of the analogues of the probability measures \eqref{eq:mun atom measure} to the unit volume Lebesgue measure
on the sphere.
On the other hand, the subtlety of the $3$-dimensional case is that, to pursue the nodal hypersurface volume, one encounters the correlations of the lattice points $\Z^{3}$ (as opposed to $\Z^{2}$), the difficulty being that
the $\Z^{3}$ does not conveniently correspond to a ring of integers or a number field, like $\Z^{2}\cong \Z[i]$. Most specifically, if one wants to analyse the nodal volume variance, then
the analogue result to Zygmund's trick \eqref{eq:Zygmund} on the length-$4$ correlations, and analogue on \eqref{eq:Rc6=o(N^4)} on the length-$6$ spectral correlations are required.
Both of these $3$d analogues were successfully resolved in ~\cite{BeMa}, by special methods, and, as a result, an analogue of \eqref{eq:var ARW KKW} was established. Building on their results on the correlations of lattice points, an analogue non-central limit theorem to \eqref{eq:Lc(n)->Meta} was established ~\cite{CamTAMS} in this case (again, converging to a single limit for the full sequence of energies).

\vspace{2mm}

The defect distribution of the ARW and deterministic toral eigenfunctions was addressed in ~\cite{KuWiYe}.
As for the spherical harmonics (\S\ref{sec:crip pnts indep intersect}), all the mentioned results are applicable for the intersection of $2$ independent ARW in $2$d ~\cite{PucattiEtAl3},
or $2$ or $3$ independent ARW in $3$d ~\cite{Notarnicola}, and some other variations of the basic problem, e.g. \cite{CamMarRosEP}. However, for dimensions higher than $3$, the results on spectral correlations are not known (and, in some cases, decisively fail), hence, save for an upper bound ~\cite{CheLa}, the problem of fluctuations of the nodal volume in dimensions $\ge 4$ is open.

\subsection{Nodal intersections against smooth curves and approximate Kac-Rice}

\label{sec:ARW ni}

Let $\Cc$ be a simple smooth curve, and $\gamma:[0,L]\rightarrow\Tb^{2}$ its unit speed parametrization, where we are interested in the number of nodal intersections $\Zc_{\Cc;n}=\Zc_{\Cc}(f_{n})$ of the ARW $f_{n}$ against $\Cc$, i.e. the number of zeros of $$f_{n}\circ \gamma = f_{n}(\gamma(\cdot)),$$ $t\in [0,L]$. Note that, unless $\Cc$ is a straight line segment, $f_{n}\circ \gamma$ is {\em not stationary}. However, since it is univariate regardless, it is easy to establish the expected number of the nodal intersections to be universally given by ~\cite{RuWiAJM}
\begin{equation*}
\E[\Zc_{\Cc;n}] = \sqrt{2n}\cdot L.
\end{equation*}
The question of the variance of $\Zc_{\Cc;n}$ is far more subtle, and, in fact, depends {\em both} on the geometry of $\Cc$ and the angular distribution of the lattice points $\Lambda_{n}$. First, it turns out that one needs to separate the two most interesting, extreme, scenarios: $\Cc$ is a straight segment, and $\Cc$ has nowhere vanishing curvature, whereas the intermediate cases could be analysed by combining both techniques.
The former case of $\Cc$ straight segment has been dealt with by Maffucci ~\cite{Maffucci}, whence various upper bounds were asserted, both conditionally and unconditionally, but no precise estimates for the variance was given for a single case.

For the other extreme case, when $\Cc$ is ``generic" simple smooth curve with {\em nowhere vanishing curvature}, one may evaluate the variance of $\Zc_{\gamma;n}$, asymptotically precisely, as follows. Let $B_{\Cc}(\Lambda_{n})$ be the quantity
\begin{equation*}
\begin{split}
B_{\Cc}(\Lambda_{n}) &:= \iint\limits_{\Cc\times\Cc} \frac{1}{N_{n}}\sum\limits_{\lambda\in\Lambda_{n}} \left\langle \frac{\lambda}{\sqrt{n}},\dot{\gamma}(t_{1})\right\rangle ^{2}
\left\langle \frac{\lambda}{\sqrt{n}},\dot{\gamma}(t_{2})\right\rangle ^{2}dt_{1}dt_{2}.
\\&=\int\limits_{0}^{L}\int\limits_{0}^{L}\int\limits_{\Sc^{1}}\langle \theta,\dot{\gamma}(t_{1}) \rangle \cdot \langle \theta,\dot{\gamma}(t_{2}) \rangle d\mu_{n}(\theta)dt_{1}dt_{2},
\end{split}
\end{equation*}
with $\mu_{n}$ as in \eqref{eq:mun atom measure}, easily generalized for an arbitrary probability measure on $\Sc^{1}$.
Observe that the expression $$(4B_{\Cc}(\Lambda_{n})-L^{2})$$ is in general fluctuating, depending on both $\gamma$ and the angular distribution of $\Lambda$, and is easy to show that for every $\Cc$, $\Lambda_{n}$ sufficiently symmetric set (equivalently, probability measure on $\Sc^{1}$),
\begin{equation*}
\frac{1}{4}L^{2} \le B_{\Cc}(\Lambda_{n}) \le \frac{1}{2}L^{2}.
\end{equation*}

Under the above notation, and assuming that $\Cc$ has nowhere vanishing curvature, it was shown ~\cite[Theorem 1.2]{RuWiAJM}, that
\begin{equation}
\label{eq:var nod intr ARW}
\var(\Zc_{\Cc;n}) = (4B_{\Cc}(\Lambda_{n})-L^{2})\cdot \frac{n}{N_{n}} + O\left(  \frac{n}{N_{n}^{3/2}} \right).
\end{equation}
In light of the above, we have $$(4B_{\Cc}(\Lambda_{n})-L^{2})\ge 0$$ always, and it does not ever vanish for a ``generic" $\Cc$. However, the vanishing of $(4B_{\Cc}(\Lambda_{n})-L^{2})$ genuinely occurs,
in two scenarios of a very different nature (still, under the nowhere vanishing curvature assumption) ~\cite[Proposition 7.1, Corollary 7.2]{RuWiAJM}: 1. Given a curve $\Cc$, $(4B_{\Cc}(\Lambda_{n})-L^{2})$ vanishes identically for every symmetric set $\Lambda_{n}$ (in a more generalized sense, for every probability measure arising from a symmetric set). In this case $\Cc$ is called ``static". An example of a static curve is a circle or a semi-circle. 2. Given a curve $\Cc$, $(4B_{\Cc}(\Lambda_{n})-L^{2})$ vanishes for precisely either the ``Cilleruelo" measure $$\frac{1}{4}\left( \delta_{\pm 1}+\delta_{\pm i}\right)$$ or its tilt by $\pi/4$ {\em tilted Cillerulo}. Otherwise, if both (1) or (2) do not hold, then the number $(4B_{\Cc}(\Lambda_{n})-L^{2})$ is bounded away from $0$.

Rossi and Wigman ~\cite{RoWi} studied the distribution of $\Zc_{\Cc;n}$ in more details, via the Wiener chaos expansion. They proved the CLT in a generic scenario: If, for some sequence $\{n\}\subseteq S$, the numbers $4B_{\Cc}(\Lambda_{n})-L^{2}$ are {\em bounded away from $0$}, then ~\cite[Theorem 1.1]{RoWi} $$\frac{\var(\Zc_{\Cc;n})-\E[\Zc_{\Cc;n}]}{\sqrt{\var(\Zc_{\Cc;n})}}   $$ converges, in distribution, to the standard Gaussian. For instance, the CLT applies in scenario (2) above, so long as one excludes the Cilleruelo and tilted Cilleruelo sequences. In this case the bulk of the fluctuations of $\Zc_{\Cc;n}$ ``lives" in the projection onto the $2$nd Wiener chaos space, with no Berry's cancellation.

On the other hand, the static curves observe a very different behaviour, both in terms of the asymptotic variance and the limit theorem. They proved ~\cite[Theorem 1.3]{RoWi} the following results, assuming that $\Cc$ is static, and under an extra assumption on the sequence $\{n\}\subseteq S$, that is generic.
Namely, it was shown that ~\cite[Theorem 1.3(1)]{RoWi}, if $\Cc$ is static, the variance of $\Zc_{\Cc;n}$ is of order of magnitude $\frac{n}{N^{2}} $, with the leading constant depending on both the angular distribution of $\Lambda_{n}$ and the geometry of $\Cc$. Further, an explicit non-Central Limit Theorem ~\cite[Theorem 1.3(2)]{RoWi} was asserted for $$\frac{\var(\Zc_{\Cc;n})-\E[\Zc_{\Cc;n}]}{\sqrt{\var(\Zc_{\Cc;n})}}   .$$ In this case Berry's cancellation occurs: The bulk of the fluctuations of $\Zc_{\Cc;n})$ ``lives" in the projection onto the $4$nd Wiener chaos space, with the projection onto the $2$nd Wiener chaos space vanishing precisely (by construction). Surprisingly enough, no analogue of static curves exist for this problem for either the random spherical harmonics or the RWM, see \S\ref{sec:nod intersect RWM}.
Finally, the analogue $3$d problems were studied: Nodal intersections against curves ~\cite{RuWiYe}, straight segments ~\cite{Maffucci2}, and surfaces ~\cite{MaRo}.

\vspace{2mm}

Note that, depending on whether $n$ is odd or even, $f_{n}$ is either odd w.r.t. to the involution $$\tau_{1}:\cdot\mapsto \cdot+(1/2,1/2)$$ or w.r.t. the involution $\tau_{2}\cdot\mapsto \cdot+(1/2,0)$ respectively.
Hence if there are sub-curves of $\Cc$ whose image, under $\tau_{1}$ or $\tau_{2}$, is lying in $\Cc$, then the sufficient conditions (i.e. non-degeneracy of the Gaussian distribution of
$(f_{n}(\gamma(t_{1})),f_{n}(\gamma(t_{2})))$, $t_{1}\ne t_{2}$) for the applicability of the Kac-Rice formula \eqref{eq:KC 2nd mom fact} fail.
We do not know whether or not such smooth curves with nowhere vanishing curvature (i.e. curves of this type that are partially invariant w.r.t. either of the $\tau_{i}$) actually exist, but, in this light, we were not able to validate whether the Kac-Rice formula \eqref{eq:KC 2nd mom fact} holds precisely in this case, i.e. whether the r.h.s. of \eqref{eq:KC 2nd mom fact} computes the $2$nd factorial moment of $\Zc_{\Cc;n}$.
However, that it does give a good
{\em estimate} for the variance, up to an admissible error term, follows from the Approximate Kac-Rice formula ~\cite[Proposition 1.3]{RuWiAJM}, especially developed for these needs, whose implementation, after much effort analysing the integral on the r.h.s. of \eqref{eq:KC 2nd mom fact}, eventually yields the asymptotic variance \eqref{eq:var nod intr ARW}. This approach also interprets the factorial (i.e. linear) part in \eqref{eq:KC 2nd mom fact} as the length of the projection onto the $1$st component of the {\em diagonal}
$$D:=\{(t,t):\: t\in [0,L]\}\subseteq [0,L]^{2},$$ corresponding to those Gaussian vectors $(f_{n}(\gamma(t_{1})),f_{n}(\gamma(t_{2})))$ that fail the non-degeneracy condition, and also gives an answer to what the r.h.s. of
\eqref{eq:KC 2nd mom fact} actually computes in case that there are degeneracies outside of $D$ that occur\footnote{Such situation will occur, for example, for $\Cc$ smooth with nowhere vanishing curvature, partially invariant w.r.t. either of $\tau_{i}$, when the answer is given in terms of a special ``degeneracy index" of $\Cc$. The author was informed that this was already known to Pavel Bleher in a different context.}.

\end{document}